\documentclass[12pt,twoside]{preprint}

\title{A new Weyl group action related to the quasi-classical Gelfand-Graev action}
\author{Xiangsheng Wang}
\institute{School of Mathematics, Shandong University,\hfill \email{wangxs1989@gmail.com}
  \\and BICMR, Peking University\hspace{3.5em}}
\date{\today}

\usepackage[english]{babel}
\usepackage[stretch=10,final]{microtype}
\usepackage[a4paper,left=2.5cm,right=2.5cm,top=4cm,bottom=2.5cm,marginparwidth=2cm]{geometry}
\usepackage{color}
\usepackage{titlesec}
\usepackage{enumitem}
\usepackage[olditem,oldenum]{paralist}
\usepackage{perpage}
\usepackage[super]{nth}
\usepackage{subdepth}
\usepackage{todonotes}
\usepackage{fancyhdr}
\usepackage[title]{appendix}
\usepackage[colorlinks,linkcolor=blue,bookmarksnumbered,pagebackref,final]{hyperref}
\usepackage{subcaption}
\usepackage[]{lipsum}
\usepackage[alphabetic,initials,msc-links]{amsrefs}

\usepackage{mathtools}
\usepackage{amsthm}
\usepackage{amssymb}
\usepackage{stmaryrd} \usepackage{mleftright}
\usepackage{tikz-cd}
\usepackage{accents}
\usepackage{nicematrix} 

\usepackage[osf,theoremfont]{newtxtext} \usepackage[upint,smallerops]{newtxmath}
\usepackage[cal=euler,scr=rsfso,frak=euler]{mathalfa}
\usepackage[T1]{fontenc}

\usepackage[normalem]{ulem}

\MakePerPage{footnote}

\DeclareSymbolFont{sfoperators}{OT1}{cmss}{m}{n}
\DeclareSymbolFontAlphabet{\mathsf}{sfoperators}
\makeatletter
\renewcommand{\operator@font}{\mathgroup\symsfoperators}
\makeatother

\numberwithin{equation}{section}
\swapnumbers                            \newtheoremstyle{plain}
{2ex plus 1ex minus .2ex}   {\medskipamount}   {\slshape}  {}       {\indent\bfseries\tlfstyle} {.}         {5pt plus 1pt minus 1pt} {}          

\newtheorem{theorem}[subsubsection]{Theorem}
\newtheorem{lemma}[subsubsection]{Lemma}
\newtheorem{corollary}[subsubsection]{Corollary}
\newtheorem{proposition}[subsubsection]{Proposition}
\newtheorem{propdef}[subsubsection]{Proposition-Definition}
\newtheorem*{claim*}{Claim}
\newtheorem*{lemma*}{Lemma}

\newtheoremstyle{definition}
{2ex plus 1ex minus .2ex}   {\medskipamount}   {}  {}       {\indent\bfseries\tlfstyle} {.}         {5pt plus 1pt minus 1pt} {}          \theoremstyle{definition}
\newtheorem{definition}[subsubsection]{Definition}
\newtheorem{remark}[subsubsection]{Remark}
\newtheorem*{remark*}{Remark}
\newtheorem*{assumption*}{Assumption}
\newtheorem*{example*}{Example}

\makeatletter
\let\newtitle\@title
\let\newauthor\@author
\let\newdate\@date
\makeatother

\titleformat{\section}{\normalfont\large\bfseries\tlfstyle}{\thesection}{0.6em}{}
\titleformat{\subsection}[runin]{\normalfont\normalsize\bfseries\tlfstyle}{\thesubsection}{0.4em}{}
\titleformat{\subsubsection}[runin]{\normalfont\normalsize\bfseries\tlfstyle}{\indent\thesubsubsection}{0.4em}{}

\usetikzlibrary{shapes,arrows}
\tikzstyle{block} = [rectangle, draw, fill=blue!20, 
    text width=5em, text centered, rounded corners, minimum height=2em]
\tikzstyle{line} = [draw, -latex']
\tikzstyle{cloud} = [draw, ellipse,fill=red!20, node distance=3cm,
    minimum height=2em]
\tikzcdset{
  arrow style=tikz,
  diagrams={>={Straight Barb[scale=0.8]}}
}

\makeatletter
\newcommand*{\transpose}{{\mathpalette\@transpose{}}}
\newcommand*{\@transpose}[2]{\raisebox{\depth}{$\m@th#1\intercal$}}
\makeatother
\newcommand{\tsp}[1]{#1^{\transpose}}

\newcommand{\Hom}[2]{\operatorname{Hom}(#1,#2)}
\newcommand{\Tr}[1]{\operatorname{Tr}(#1)}
\newcommand{\Spec}[1]{\operatorname{Spec}(#1)}
\newcommand{\aset}[1]{\{#1\}}
\newcommand{\ubar}[1]{\underaccent{\bar}{#1}}
\newcommand{\Ad}[2]{\operatorname{Ad}_{#1}(#2)}
\newcommand{\diag}[1]{\operatorname{diag}(#1)}
\newcommand{\rGL}{\mathrm{GL}}
\newcommand{\rSL}{\mathrm{SL}}
\newcommand{\rU}{\mathrm{U}}
\newcommand{\rSU}{\mathrm{SU}}
\newcommand{\kgl}{\mathfrak{gl}}
\newcommand{\ksl}{\mathfrak{sl}}
\newcommand{\longeq}{=\joinrel=\joinrel=}
\newcommand{\simarrow}{\xrightarrow{
    \smash{\raisebox{-0.65ex}{\ensuremath{\scriptstyle\sim}}}}}
\newcommand{\gr}{\mathrm{gr}}
\newcommand{\zm}[2]{\mathbf{0}_{{#1},{#2}}}

  \newcommand{\ka}{\mathfrak{a}}
\newcommand{\kb}{\mathfrak{b}}
\newcommand{\kd}{\mathfrak{d}}

\newcommand{\kg}{\mathfrak{g}}
\newcommand{\kt}{\mathfrak{t}}

\newcommand{\ku}{\mathfrak{u}}
\newcommand{\kM}{\mathfrak{M}}
\newcommand{\cB}{\mathcal{B}}

\newcommand{\cM}{\mathcal{M}}
\newcommand{\cZ}{\mathcal{Z}}

\newcommand{\rA}{\mathrm{A}}

\newcommand{\rT}{\mathrm{T}}
\newcommand{\rM}{\mathrm{M}}

\newcommand{\rr}{\mathrm{r}}

\newcommand{\rv}{\mathrm{v}}
\newcommand{\opi}{\mathrm{i}}

\newcommand{\rsc}{\mathrm{sc}}
\newcommand{\bZ}{\mathbb{Z}}
\newcommand{\bC}{\mathbb{C}}
\newcommand{\bR}{\mathbb{R}}

\newcommand{\bW}{\mathbb{W}}
\newcommand{\bO}{\mathbb{O}}
\newcommand{\bfv}{\mathbf{v}}
\newcommand{\bfw}{\mathbf{w}}

\newcommand{\bfA}{\mathbf{A}}
\newcommand{\bfC}{\mathbf{C}}
\newcommand{\bfI}{\mathbf{I}}
\newcommand{\bfM}{\mathbf{M}}
\newcommand{\bfN}{\mathbf{N}}
\newcommand{\sfs}{\mathsf{s}}

\newcommand{\sch}{\substack{h\in Q \\ \hd{h} = k}}
\newcommand{\sct}{\substack{h\in Q \\ \tl{h} = k}}
\newcommand{\ade}{\mathrm{ADE}}
\newcommand{\aff}{\mathrm{aff}}
\newcommand{\grk}{{\rmfamily G\textsuperscript{3}RK}}
\newcommand{\afb}{(G/U)_{\aff}}
\newcommand{\afbc}{\big(\rT^*(G/U)\big)_{\aff}}
\newcommand{\DQ}{\mathrm{D}Q}
\newcommand{\DQw}{\DQ^{\bfw}}
\newcommand{\pt}{(3)}
\newcommand{\ras}{\mathrm{rs}}
\newcommand{\bkb}{\bar{\mathfrak{b}}}
\newcommand{\ukb}{\Ad{u^{-1}}{\bkb}}

\newcommand{\dsl}{\mathbin{/\mkern-5mu/}}

\newcommand{\kah}{K\"ahler}
\newcommand{\hkh}{hyper-K\"ahler}
\newcommand{\SL}[1]{\mathrm{SL}(#1,\mathbb{C})}

\newcommand{\U}[1]{\mathrm{U}(#1)}
\newcommand{\GL}[1]{\mathrm{GL}(#1,\mathbb{C})}
\newcommand{\PGL}[1]{\mathrm{PGL}(#1,\mathbb{C})}
\newcommand{\hd}[1]{\operatorname{h}(#1)}
\newcommand{\tl}[1]{\operatorname{t}(#1)}
\newcommand{\op}{\mathrm{op}}
\newcommand{\surj}{\mathrm{surj}}

\DeclareMathOperator{\diff}{d}

\DeclareMathOperator{\opj}{j}

\DeclareMathOperator{\id}{id}
\DeclareMathOperator{\im}{in}
\DeclareMathOperator{\om}{out}
\DeclareMathOperator{\opp}{p}

\DeclareMathOperator{\opN}{N}

\DeclareMathOperator{\mm}{m}

\begin{document}

\maketitle
\begin{abstract}
  We construct a Weyl group action on the DKS type varieties, a certain class of varieties associated with quivers.
  As a result, on some special DKS type varieties, we can give a quiver theoretic explanation of the quasi-classical Gelfand-Graev action discovered by Ginzburg and Riche and studied by Ginzburg and Kazhdan recently.
\end{abstract}

\section{Introduction}

\subsection{Backgrounds.}
Let $G$ be a complex connected semi-simple group and $U$ be a maximal unipotent subgroup.
The homogeneous manifold $G/U$ is sometimes called the ``basic (base) affine space'' of the group $G$, c.f.~\cite{Bernv-steui-n:1975bp,Bezrukavnikov:2002aa}.
The importance of this space is partially reflected by the following empirical fact summarized in~\cite{Bernv-steui-n:1975bp}:
\begin{quote}
  \emph{``Experience in representation theory suggests that for many problems in representation theory the solution results from a careful study of the base affine space.''}
\end{quote}
Guided by this philosophy, in this paper, we are interested in investigating a natural group action on the affinization of $\rT^* (G/U)$.

By studying the ring of algebraic differential operators on $G/U$, in an unpublished paper written in the 1960s, S.\ Gelfand
and M.\ Graev discover a nontrivial action of the Weyl group on this ring.
For example, when the Weyl group is $\bZ_2$, such an action is essentially the Fourier transformation on the polynomial differential operators.
In view of this example, it seems that one can hardly expect such an action to be induced from a Weyl group action on $G/U$ itself.
Therefore, it is understandable that Ginzburg and Riche~\cite{Ginzburg:2015so} describe it as ``quite surprising'' when they find a natural Weyl group action on $\afbc$, the affinization of the cotangent bundle of $G/U$.
Based on this new action, called the \grk\ action in the following, Ginzburg and Kazhdan~\cite{Ginzburg:2018la} give an algebraic interpretation of the Gelfand-Graev action.
Due to this fact, a reasonable way to view the \grk\ action is treating it as a quasi-classical analog of the Gelfand-Graev action.
In fact, in~\cite{Ginzburg:2018la}, the \grk\ action is called the \emph{quasi-classical Gelfand-Graev action} directly.
Moreover, by using the canonical map from $\rT^*(G/U)$ to $\tilde{\kg}$, the Grothendieck-Springer resolution of $\kg$, the \grk\ action is also related to another famous Weyl group action: the Springer representation~\cite{Springer:1976wf,Springer:1978ye}.

Under a different circumstance, in symplectic geometry or \kah\ geometry, the basic affine space and its cotangent bundle also arise.
Let $K$ be a maximal compact subgroup of $G$.
In~\cite{Guillemin:2002sx}, Guillemin, Jeffrey and Sjamaar study a construction called symplectic implosion on Hamiltonian $K$-manifolds, which shares similar properties with the more familiar symplectic reduction construction.
To get an impression of the relation between the symplectic implosion and the basic affine space, one can perform the implosion construction on $\rT^*K \simeq G$, which gives rise to $\afb$, the affinization of the quasi-affine space $G/U$.
Similar to what happens in the representation theory, $\afb$ also occupies a special position in the implosion construction, i.e.\ as a certain ``universal'' object.
Moreover, when $G= \SL{n+1}$, it is possible to generalize the implosion construction to \hkh\ manifolds,~\cite{Dancer:2013aa,Dancer:2016ab}.
In particular, in such a generalization, the space $\afbc$ is the counterpart of $\afb$.

An interesting point is that if we study $\afbc$ using the symplectic geometry method, this space is also related to some Weyl group action, however in an indirect way.
To see this, we recall that Dancer, Kirwan and Swann undertake a detailed analysis of $\afbc$ by realizing it as a variety associated with the $\rA_{n}$ type quiver, c.f.~\cite[Theorem~7.18]{Dancer:2013aa}.
Then, by using the quiver theoretic structure on $\afbc$, the \hkh\ quotient of the Hamiltonian $(\U{1})^{n}$-action on $\afbc$ is the Nakajima quiver variety $\kM_{\zeta}(\bfv,\bfw)$ exactly, where $\bfv = \tsp{(n,n-1,\cdots,1)}$, $\bfw = \tsp{(n+1,0,\cdots,0)}$ are dimension vectors.
Now, by a result of Nakajima,~\cite[Proposition~9.1]{Nakajima:1994aa}, as we recall in Theorem~\ref{thm:nweyl}, the Weyl group has a natural action on the set $\{\kM_{\zeta}(\bfv,\bfw) | \zeta\in \bR^3\otimes \opi \bR^{n}\}$.

\subsection{Main results.}
\label{sub:comp-q}
Comparing the above two scenarios that $\afbc$ appears, the following question occurs to us quite naturally.
Since there are two Weyl group actions associated with $\afbc$ in some way, one by the \grk\ action, another via its quiver varieties explanation and Nakajima's result, do they correlate with each other?
A classical result pertains to this question is that the Springer representation for the Weyl group has an interpretation in terms of the Weyl group action on the Nakajima quiver varieties,~\cites{Nakajima:1994aa,Slodowy:1980aa}.
In view of this fact, as well as the relation between the Springer representation and the \grk\ action, the answer to the above question seems to be positive.

Note that Nakajima's result only deals with the \hkh\ quotients of $\afbc$.
To give a precise formulation of our comparison question, we would like to construct a Weyl group action on $\afbc$ itself, not on its \hkh\ quotients, in a quiver theoretic way.
However, since $\afbc$ is just an example of a certain class of varieties associated with a quiver, which are called DKS type varieties, if we can find a quiver theoretic Weyl group action on $\afbc$, it seems equally plausible to expect such an action existing on a general DKS type variety.
This is true indeed.
Given the dimension vectors $\bfv,\bfw$ of a quiver, we denote such a variety by $\cM(\bfv,\bfw)$.
We are going to show the following result.

\begin{theorem}[{$=$Theorem~\ref{thm:main-re}}]
  \label{thm:main}
  Fix a quiver $Q$ of $\ade$ type, i.e.\ having a finite Weyl group.
  Let $\bfC$ be the Cartan matrix of the quiver.
  If $\bfw = \bfC \bfv$, then there is a Weyl group action on the variety $\cM(\bfv,\bfw)$ associated with $Q$.
\end{theorem}

The Weyl group action in the above theorem is compatible with the Weyl group action on $\{\kM_{\zeta}(\bfv,\bfw) | \zeta\in \bR^3\otimes \opi \bR^{n}\}$, i.e.\ Nakajima's theorem mentioned before.
In fact, note that the \hkh\ reduction of $\cM(\bfv,\bfw)$ at $\zeta$ is $\kM_{\zeta}(\bfv,\bfw)$.
Via such a reduction, the Weyl group action on $\cM(\bfv,\bfw)$ descends exactly to the Weyl group action appearing in Nakajima's theorem.
Therefore, the Weyl group action in the above theorem fits well for our comparison question raised at the beginning of this subsection.

Our method to show Theorem~\ref{thm:main} is based on the algebraic approach to proving Nakajima's theorem, i.e.\ by using the reflection functor, c.f.~\cites{Lusztig:2000aa,Maffei:2002aa,Nakajima:2003aa}.
However, the reflection functor in the literature is only devised for the usual quiver varieties, which, roughly speaking, are some quotient spaces associated with the general linear groups.
As a contrast, the varieties that we consider in this paper, $\cM(\bfv,\bfw)$, are only some quotient spaces associated with the special linear groups.
Therefore, we need to check the well-known method to construct the reflection functor also works for $\cM(\bfv,\bfw)$.
More precisely, we show that even with additional constraint coming from the torsion of chain complexes, it is still possible to construct the reflection functor.

Once we establish the existence of the Weyl group action on $\cM(\bfv,\bfw)$, we can compare it with the \grk\ action and answer the comparison question affirmatively for some special varieties.
Here is a result of this type.

\begin{theorem}[{$=$Theorem~\ref{thm:main2-re}}]
  \label{thm:main2}
  Choose the dimension vectors to be $\bfv = \tsp{(n,n-1,\cdots,1)}$, $\bfw = \tsp{(n+1,0,\cdots,0)}$ for a quiver of $\mathrm{A}_n$ type.
  Take $G = \SL{n+1}$ and let $\Xi: \cM{(\bfv,\bfw)} \rightarrow \afbc$ be the isomorphism constructed in~\cite{Dancer:2013aa}.
Then $\Xi$ induces an equivariant map with respect to the Weyl group action on $\cM(\bfv,\bfw)$ and the \grk\ action on $\big(\rT^*(G^{\mathrm{ad}}/U)\big)_{\aff}$, where $G^{\mathrm{ad}} = \PGL{n+1}$ is the adjoint group associated with $G$.
\end{theorem}

We remark that it is also true that the isomorphism $\Xi$ itself is equivariant.
To show this result, we need to deal with more technicalities.
But the overall idea behind the proof is the same with Theorem~\ref{thm:main2}.
Therefore, to make the argument more transparent, we choose to provide all the details for the proof of Theorem~\ref{thm:main2} solely.
After that, in Remark~\ref{rk:pf-sc}, we discuss how to modify the proof of Theorem~\ref{thm:main2} to show that $\Xi$ is equivariant.

\subsection{An outline of the paper.}
This paper is organized as follows.
In Section~\ref{sec:quiver-varieties}, we review some backgrounds about quiver varieties and set up notations used in the whole paper.
In Section~\ref{sec:dks-type-varieties}, we define the DKS type variety and show that it carries an action of the Weyl group.
Especially, we restate Theorem~\ref{thm:main} in Theorem~\ref{thm:main-re} and prove Theorem~\ref{thm:main-re}.
In Section~\ref{sec:comparison}, we recall the \grk\ action on $\afbc$.
For $G = \SL{n+1}$, we review a description of $\afbc$ using the DKS type varieties after the method in~\cite{Dancer:2013aa}.
Then, we restate Theorem~\ref{thm:main2} in Theorem~\ref{thm:main2-re} and prove Theorem~\ref{thm:main2-re}.

\subsection*{Acknowledgement.}
\label{sec:ack}
The author has been supported by China Postdoctoral Science Foundation (Grant No.\ BX201700008), and the fundamental research funds of Shandong University (Grant No.\ 2020GN063).
The author is very grateful for the helpful correspondence with Prof.\ Victor Ginzburg about the contents of this paper and he also would like to thank Prof.\ Gang Tian and Prof.\ Weiping Zhang for their encouragement during the preparation of this paper.
Last but not least, the author thanks the anonymous referee, whose comments help to improve the quality of this paper to a great degree.

\section{Quiver varieties}\label{sec:quiver-varieties}
\subsection{}
In this section, we review some materials about quiver varieties and set up notations used in this paper.
A useful resource for the material covered here is Ginzburg's lectures~\cite{Ginzburg:2012aa} on the Nakajima quiver varieties.

We are always working over the field $\bC$.
By a \emph{variety} we mean that it is a separated scheme of finite type over $\bC$.
And we only consider the closed points of a scheme.

Given a symplectic manifold $(M,\omega)$ with a Hamiltonian group action $G$, we define the (real) moment map $\mu: M \rightarrow \kg^*$ of this group action \emph{with the sign convention} such that
\begin{equation*}
  \iota_{V^M} \omega = \diff \langle \mu, V \rangle,
\end{equation*}
where $V\in \kg$ and $V^M$ is the vector field generated by $V$.
The complex moment map is defined using the similar sign convention.

\subsection{A review of quiver varieties.}
\label{sub:quiver}
By a \emph{quiver} $Q$, we mean that it is a finite directed graph with the vertex set $I$ and the edge set $E$.
In the following, we will identify $I$ with the set $\{1,\cdots,n\}$, where $n = \#I$.
Each element $h\in E$ has an orientation, i.e.\ $h$ is an arrow.
The head vertex and the tail vertex of $h$ are denoted by $\hd{h}$ and $\tl{h}$ respectively.
We will always assume that $Q$ has no edge loops, i.e.\ no edges joining a vertex with itself.
The \emph{opposite quiver} $Q^{\op}=(I,\bar{E})$ of $Q$ is isomorphic to $Q$ as an undirected graph, but the orientation of each edge is reversed.
For $h\in E$, we denote $\bar{h}\in \bar{E}$ to be the edge corresponding to $h$ but with reversed orientation.
The \emph{double quiver} of $Q$ is defined by $\DQ \coloneqq (I,E \sqcup \bar{E})$.
To save symbols, if there is no ambiguity, for a quiver $Q = (I, E)$, $h\in E$ is often written as $h\in Q$.
\begin{figure}[htbp]
  \centering
  \begin{equation*}
    \begin{tikzcd}
     \phantom{a} & \bullet \arrow[l, phantom, "{\scriptstyle{\tl{h}}}" near start] \arrow[r, "h"] & \bullet \arrow[r, phantom, "{\;\;\scriptstyle{\hd{h}}}" near start] & \phantom{a} \\
     \phantom{a} & \bullet \arrow[l, phantom, "{\scriptstyle{\hd{\bar{h}}}\;\;}" near start] & \arrow[l,"\bar{h}"] \bullet \arrow[r, phantom, "{\scriptstyle{\tl{\bar{h}}}}" near start] & \phantom{a}
    \end{tikzcd}
  \end{equation*}
  \caption{Two orientations for an edge.}
\end{figure}
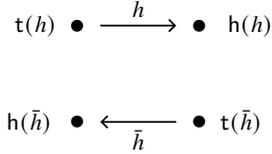

Let $\bfA = (\bfA_{kl})_{k,l\in I}$ be the \emph{adjacency matrix} of $\DQ$, where $\bfA_{kl}$ is the number of edges of $\DQ$ going from $l$ to $k$.
Clearly, $\bfA$ is also the adjacency matrix for the underlying undirected graph of $Q$.
The symmetric \emph{generalized Cartan matrix} associated with $Q$ is defined as $\bfC \coloneqq 2\bfI - \bfA$, where $\bfI$ is the identity matrix.

Using $\bfC$, one can construct a Kac-Moody algebra, whose Weyl group $\bW$ is a discrete group defined by generators $\{s_k\}_{k\in I}$ satisfying the following relations,
\begin{subequations}
  \begin{alignat}{2}
    & s_k^2 = 1, &&  \label{eq:cox1}\\
    & s_ks_l = s_ls_k, &&\quad \text{if } \bfA_{kl} = 0, \label{eq:cox2}\\
    & s_ks_ls_k = s_ls_ks_l, &&\quad \text{if } \bfA_{kl} = 1. \label{eq:cox3}
  \end{alignat}
\end{subequations}
Note that the group $\bW$ is finite if and only if the underlying undirected graph of $Q$ is a Dynkin diagram of $\ade$ type.

On $\bZ^n$, we can define a $\bW$ action as follows.
Taking $\zeta = \tsp{(\zeta_1,\cdots,\zeta_n)} \in  \bZ^n$,
\begin{equation}
  \label{eq:weylac}
  s_k(\zeta_1,\cdots,\zeta_l,\cdots,\zeta_n) \coloneqq (\zeta_1 - \bfC_{k1}\zeta_k,\cdots,\zeta_l - \bfC_{kl}\zeta_k,\cdots,\zeta_n - \bfC_{kn}\zeta_k).
\end{equation}
For the later usage, we also extend the above action to an action on $\bR^n$ or $\bC^n$.

To each vertex $k$, we associate two $\bC$ vector spaces $V_k,W_k$.
Let $V = (V_k)_{k\in I}$ and $W = (W_k)_{k\in I}$.
The following column vectors are called dimension vectors of $V$, $W$ respectively,
\begin{align*}
  \bfv &\coloneqq \tsp{(\dim{V_1}, \dim{V_2},\cdots, \dim{V_n})}\in \bZ^n_{\ge 0},\\
  \bfw &\coloneqq \tsp{(\dim{W_1}, \dim{W_2},\cdots, \dim{W_n})} \in \bZ^n_{\ge 0}.
\end{align*}
In the whole paper, the symbol $\dim$ denotes the complex dimension of a complex vector space.
And we often use $v_k$ (resp.\ $w_k$) as an abbreviation for $\dim{V_k}$ (resp.\ $\dim{W_k}$).
Following Nakajima,~\cite{Nakajima:1994aa}, we define the following vector space,
\begin{equation}
  \label{eq:defm}
  \bfM(V,W) \coloneqq \Big(\bigoplus_{h\in \DQ} \Hom{V_{\tl{h}}}{V_{\hd{h}}}\Big) \oplus \Big(\bigoplus_{k\in I} \Hom{W_k}{V_k} \oplus \Hom{V_k}{W_k} \Big).
\end{equation}
When we are only concerned with the isomorphism classes of $V,W$, we use the notation $\bfM(\bfv,\bfw)$, and even $\bfM$ for the same space in (\ref{eq:defm}).
Besides, for the $\bfw = 0$ case that we consider at times, we also use the notation $\bfM(\bfv)$ for $\bfM(\bfv,0)$.
An element of $\bfM(V,W)$ is usually denoted by $(B,i,j)$, where $B = (B_h)_{h\in \DQ}$, $i = (i_k)_{k\in I}$ and $j = (j_k)_{k\in I}$.

We can naturally view $\bfM(V,W)$ as a cotangent bundle over a complex space, c.f.~\cite[(4.3.3)]{Ginzburg:2012aa},~\cite[(2.4)]{Nakajima:1994aa}.
Therefore, there is a canonical (complex) symplectic form on $\bfM(V,W)$,
\begin{equation}
  \label{eq:sympformm}
  \omega \big((B,i,j),(B',i',j')\big) = \sum_{h\in Q}\big(\Tr{B_hB'_{\bar{h}}} - \Tr{B_{\bar{h}}B'_{h}}\big) + \sum_{k\in I} \big(\Tr{j_ki'_k} - \Tr{i_kj'_k}\big).
\end{equation}
With respect to $\omega$, the group $\rGL_{\bfv} \coloneqq \prod_{k\in I} \mathrm{GL}(V_k)$ acts on $\bfM$ symplectically,
\begin{equation*}
  (B_h,i_k,j_k) \mapsto (g_{\hd{h}}B_hg_{\tl{h}}^{-1},g_{k}i_k,j_kg_{k}^{-1}).
\end{equation*}
Moreover, there exists a (complex) moment map for the above $\rGL_{\bfv}$ action.
We denote the Lie algebra of $\rGL_{\bfv}$ by $\kgl_{\bfv} \coloneqq \oplus_{k\in I} \kgl(V_k)$, whose dual is denoted by $\kgl_{\bfv}^*$.
By identifying $\kgl_{\bfv}$ with $\kgl_{\bfv}^*$ using the Killing form (the sum of traces of each component), the moment map is written as follows,
\begin{equation}
  \label{eq:mm}
  \mu(B,i,j) = (\sum_{\sch} B_hB_{\bar{h}} - \sum_{\sct} B_{\bar{h}}B_{h} + i_kj_k)_k \in \bigoplus_{k\in I} \kgl(V_k) = \kgl_{\bfv}.
\end{equation}

If we further equip $V_k,W_k$ with Hermitian inner products, then it is possible to define a Hamiltonian \hkh\ structure on $\bfM$ whose complex symplectic form is $\omega$ and complex moment map is $\mu$.
Since we will not use this \hkh\ structure, we refer readers to~\cite[\S~2]{Nakajima:1994aa} for details.
But as a helpful comparison with $\mu$, we write down the real moment map of the Hamiltonian \hkh\ structure with respect to the group $\rU_{\bfv} \coloneqq \prod_{k\in I} \rU(V_k) \subseteq \rGL_{\bfv}$.
Denote the Lie algebra of $\rU(V_k)$ by $\ku_{\bfv}\coloneqq \oplus_{k\in I} \ku(V_k)$ and identify $\ku_{\bfv}$ with its dual using the Killing form.
\begin{equation}
  \label{eq:rmm}
  \mu_{\bR}(B,i,j) = \frac{\opi}{2}(\sum_{\substack{h\in \DQ \\ \hd{h}=k}} B_hB_h^{\dagger} - B^{\dagger}_{\bar{h}}B_{\bar{h}} + i_ki_k^{\dagger} - j_k^{\dagger}j_k)_k\in \bigoplus_{k\in I} \ku(V_k) = \ku_{\bfv},
\end{equation}
where $(\bullet)^{\dagger}$ is the Hermitian adjoint.

Since the center of each summand in $\kgl(V_k)$ (resp.\ $\ku(V_k)$) is a 1-dimensional Lie algebra of scalar matrices, we can and will identify the center of $\kgl_{\bfv}$ (resp.\ $\ku_{\bfv}$) with $\bC^n$ (resp.\ $\bR^n$).
Taking $\zeta = (\lambda,\theta) \in \bC^n \oplus \bR^n$, the \emph{Nakajima quiver variety} is defined as the following \hkh\ quotient,
\begin{equation}
  \label{eq:n-v}
  \kM_{\zeta} = \kM_{\zeta}(\bfv,\bfw) \coloneqq \big(\mu^{-1}(\lambda) \cap \mu^{-1}_{\bR}(\theta)\big)/\rU_{\bfv}.
\end{equation}
As usual, we can also define $\kM_{\zeta}$ as a GIT quotient.
For simplicity, assuming $\theta = 0$, the result is
\begin{equation}
  \label{eq:n-vc}
  \kM_{(\lambda,0)} = \mu^{-1}(\lambda) \dsl \rGL_{\bfv} = \Spec{\bC[\mu^{-1}(\lambda)]^{\rGL_{\bfv}}}.
\end{equation}

\medskip
\noindent \uline{A clarification of notations}.
Note that we have chosen the notation $\kM_{\zeta}(\bfv,\bfw)$ for the Nakajima quiver varieties.
By our general rule for notations, it implies that $\kM_{\zeta}$ only depends on the dimension vectors up to isomorphism, which is correct but hides some delicate details.
It is clear from the definition that for $V,V'$ with the same dimension vector $\bfv$, there exists a canonical isomorphism from $\kM_{\zeta}(V,W)$ to $\kM_{\zeta}(V',W)$.
But for $W,W'$ with the same dimension vector $\bfw$, the isomorphism from $\kM_{\zeta}(V,W)$ to $\kM_{\zeta}(V,W')$ in fact depends on the choice of isomorphism between $W$ and $W'$.
In this paper, we are only concerned with isomorphism classes of $V$ actually.
We will simply fix $W = (W_k)_{k\in I}$ once and for all to avoid such ambiguities.
\medskip

In~\cite[Definition~2.3]{Nakajima:2003aa}, fixing $\bfw$, Nakajima defines the following affine action of the Weyl group $\bW$ on $\bZ^n$, the lattice spanned by dimension vectors.
Let $\bfv = \tsp{(v_1,\cdots,v_n)}, \bfw = \tsp{(w_1,\cdots,w_n)}$.
\begin{equation}
  \label{eq:affweyl}
  s_k *_{\bfw} (v_1,\cdots,v_k,\cdots,v_n) \coloneqq (v_1,\cdots, \overbrace{v_k - \textstyle\sum_{l=1}^n\bfC_{kl}v_l + w_k}^{k\text{th entry}},\cdots,v_n).
\end{equation}
Note that if $\bfw = 0$, the above action is conjugate to the action defined in (\ref{eq:weylac}).
In fact, for these two actions, we have the relation $\bfw - \bfC(s_k *_{\bfw} \bfv) = s_k(\bfw - \bfC \bfv)$.
With this affine action, Nakajima~\cite[Proposition~9.1]{Nakajima:1994aa} constructs the following collection of isometries on the quiver varieties.
\begin{theorem}
  \label{thm:nweyl}
  Let $\zeta$ be a generic point in $\bC^n \oplus \bR^n$.
  If $\bfC$ is positive definite, i.e.\ $Q$ is a finite quiver of Dynkin $\ade$ type, for any $\sigma\in \bW$, there exists a \hkh\ isometry
  \begin{equation*}
    \Phi_{\sigma}:\kM_{\zeta}(\bfv,\bfw) \rightarrow \kM_{\sigma(\zeta)}(\sigma *_{\bfw} \bfv, \bfw).
  \end{equation*}
  Moreover, if $\tau$ is another element of $\bW$, we have $\Phi_{\sigma}\Phi_{\tau} = \Phi_{\sigma\tau}$.
\end{theorem}
For other constructions of the collection of isometries in the above theorem, as well as some generalizations, we refer readers to~\cite{Crawley-Boevey:1998aa,Lusztig:2000aa,Maffei:2002aa,Nakajima:2003aa}.

We should remind readers that the affine action of $\bW$ on $\bZ^n$ in (\ref{eq:affweyl}) generally does not preserve $\bZ^n_{\ge 0}$.
In this paper, we are only interested in a special case, i.e.\ $\bfw = \bfC \bfv$.
With this condition, one can check that $\bfv = \sigma *_{\bfw} \bfv$ for any $\sigma\in \bW$.

\subsection{Crawley-Boevey's construction.}\label{sub:cb}
By considering the vector spaces $\{W_k\}_{k\in I}$, we actually introduce a larger quiver $Q^{\heartsuit} \supseteq Q$, which is sometimes called the \emph{framing} of $Q$,~\cite[\S~3]{Ginzburg:2012aa}.
If $\bfw = 0$, the Nakajima quiver varieties degenerate to the \emph{quiver varieties without framings}, which are the usual meanings of ``quiver varieties'' in the literature.
For the later usage, we recall a construction of the Nakajima quiver varieties by using the quiver varieties without framings due to Crawley-Boevey~\cite[p.~261]{Crawley-Boevey:2001aa}.

As before, let $Q = (I,E)$ be a quiver and $\bfv,\bfw$ be dimension vectors.
Crawley-Boevey considers a quiver $Q^{\bfw}$ with the vertex set $I \sqcup \{\infty\}$, where $\infty$ is a newly added vertex.
For each vertex $k\in I$, we draw $w_k$ edges from $k$ to $\infty$.
The edge set of the quiver $Q^{\bfw}$ is the union of $E$ and the newly drawn edges.
We define a new dimension vector $\hat{\bfv}\in \bZ^{I\sqcup \aset{\infty}}_{\ge 0}$ of $Q^{\bfw}$, that is, $\hat{v}_k \coloneqq v_k$ for $k\in I$ and $\hat{v}_{\infty} \coloneqq 1$.
There is a natural group inclusion $\rGL_{\bfv} \hookrightarrow \rGL_{\hat{\bfv}}$ that sends an element $g = (g_k)_{k\in I}$ to $\hat{g} = (\hat{g}_k)_{k\in I\sqcup \aset{\infty}}$, where $\hat{g}_k \coloneqq g_k$ for $k\in I$ and $g_{\infty} \coloneqq \mathrm{Id}$.
As in the $\bfM(\bfv,\bfw)$ case, $\bC^{*}$, as a subgroup of $\rGL_{\hat{\bfv}}$, can act on $\bfM(\hat{\bfv})$. However, the $\bC^{*}$ action defined in this way is trivial actually.
Therefore, via the isomorphism $\rGL_{\bfv} \simeq \rGL_{\hat{\bfv}}/\bC^{\ast}$, we can and will treat $\bfM(\hat{\bfv})$ as a $\rGL_{\bfv}$-variety.
Let $\hat{\mu}$ be the moment map of $\bfM(\hat{\bfv})$.
By using $\bfM(\hat{\bfv})$ and $\hat{\mu}$, the quiver varieties without framing associated with $Q^{\bfw}$ are defined in the same way as (\ref{eq:n-v}) or (\ref{eq:n-vc}).
Therefore, to compare the quiver varieties without framing associated with $Q^{\bfw}$ and the Nakajima quiver varieties associated with $Q$, we only need to compare the spaces ($\bfM(\hat{\bfv})$ and $\bfM(\bfv,\bfw)$) and the moment maps ($\hat{\mu}$ and $\mu$).

\begin{figure}[htbp]
  \centering
  \begin{minipage}{150pt}
    \centering
    \begin{tikzcd}
      W_1 & & W_2 \\[-2em]
      \circ & & \circ \\
      \bullet \arrow[u] \arrow[rr,dashed] & & \bullet \arrow[u]\\[-2em]
      V_1  & & V_2
    \end{tikzcd}
    \subcaption{$Q^{\heartsuit}$}
  \end{minipage}
  \begin{minipage}{150pt}
    \centering
    \begin{tikzcd}
      & \bC & \\[-2em]
      & \infty & \\
      \bullet \arrow[ur, bend left = 20, ""{name=L1}] \arrow[ur, bend right = 20, ""{name=L2,below}] \arrow[rr, dashed] & & \bullet \arrow[ul, bend left = 20, ""{name=R1}]
      \arrow[ul, bend right = 20, ""{name=R2,above}] \arrow[dashed, no head, from=L1, to=L2, "w_1"] \arrow[dashed, no head, from=R2, to=R1, "w_2"] \\[-2em]
      V_1 & & V_2
    \end{tikzcd}
    \subcaption{$Q^{\bfw}$}
  \end{minipage}
\caption{An example of the correspondence between $Q^{\heartsuit}$ and $Q^{\bfw}$.}
\label{fig:qqw}
\end{figure}

We introduce a vector space $\hat{V} \coloneqq \oplus_{I\sqcup \aset{\infty}} V_k$, where $V_{\infty}$ is identified with $\bC$ (1-dimensional vector space with a selected basis vector).
To compare $\bfM(\hat{\bfv})$ and $\bfM(\bfv,\bfw)$, we also need to pick out a basis for each vector space $W_k$ of $W$.
Now let $\aset{h_{\alpha}}_{\alpha = 1}^{w_k}$ be the set of edges satisfying $\hd{h_{\alpha}} = \infty, \tl{h_{\alpha}} = k$.
By using the chosen basis of $W_k$, we can find a bijection between two sets: $\aset{(B_{h_{\alpha}})_{\alpha = 1}^{w_k}| B_{h_{\alpha}}\in \Hom{V_k}{V_{\infty}}\simeq \Hom{V_k}{\bC}}$ and $\Hom{V_k}{W_k}$.
In the same fashion, there is a bijection between $\Hom{W_k}{V_k}$ and $\aset{(B_{h_{\alpha}})_{\alpha = 1}^{w_k}| B_{h_{\alpha}}\in \Hom{V_{\infty}}{V_k}}$.
Due to these two bijections, we have a $\rGL_{\bfv}$-equivariant isomorphism
\begin{equation*}
  \Delta:\bfM(V,W) \simarrow \bfM(\hat{V}).
\end{equation*}
The moment map $\hat{\mu}$ of $\bfM(\hat{V})$ is a little cleaner than (\ref{eq:mm}). For $\hat{B}\in \bfM(\hat{V})$,
\begin{equation}
  \label{eq:mm-v}
  \hat{\mu}(\hat{B}) = (\sum_{\substack{h\in Q^{\bfw} \\ \hd{h} = l}} \hat{B}_h\hat{B}_{\bar{h}} - \sum_{\substack{h\in Q^{\bfw} \\ \tl{h} = l}} \hat{B}_{\bar{h}}\hat{B}_{h})_l \in (\bigoplus_{k\in I} \kgl(V_k)) \oplus \bC = \kgl_{\hat{\bfv}}.
\end{equation}
Moreover, since the trace vanishes on commutators, $\hat{\mu}$ must take values in the following subspace of $\kgl_{\hat{\bfv}}$,
\begin{equation*}
  \kgl'_{\bfv}\coloneqq \aset{(L_l)_{l\in I\sqcup \aset{\infty}}\in \kgl_{\hat{\bfv}}| \sum_{l\in I\sqcup \aset{\infty}} v_l\cdot\Tr{L_l} = 0}.
\end{equation*}
Let $\opj$ be the following isomorphism,
\begin{equation*}
  \opj(L_1.\cdots,L_n) \coloneqq (L_1,\cdots,L_n, - \sum_{k\in I} v_k\cdot\Tr{L_k}): \kgl_{\bfv}\rightarrow \kgl'_{\bfv}.
\end{equation*}
The relation between $\mu$ and $\hat{\mu}$ is
\begin{equation*}
  \hat{\mu}\circ \Delta = \opj \circ \mu.
\end{equation*}
Therefore, we often identify $\kgl_{\bfv}$ and $\kgl_{\bfv}'$ and use the same notation for moment maps on $\bfM(V,W)$ and $\bfM(\hat{V})$ when there is no ambiguity.

Denote the Weyl group of $Q^{\bfw}$ by ${\bW}^{\bfw}$.
By definition, ${\bW}^{\bfw}$ contains $\bW$ as a subgroup.
For any $\bfv \in \bZ^n$, $k\in I$, one can check the following equality,
\begin{equation*}
  s_{k,{\bW}^{\bfw}}*_{0}\hat{\bfv} = \widehat{s_{k,\bW} *_{\bfw} \bfv}.
\end{equation*}
That is, the map $\bfv \mapsto \hat{\bfv}$ is equivariant with respect to the affine action of $\bW$.

In summary, by identifying some spaces, as well as some moment maps, as above, we know that the quiver varieties without framings associated with $Q^{\bfw}$ are just another construction of the Nakajima quiver varieties associated with $Q$.
Besides, when discussing the Weyl group action on the Nakajima quiver varieties, these two constructions using $Q$ or $Q^{\bfw}$ lead to equivalent results.

\section{DKS type varieties and their Weyl group action}\label{sec:dks-type-varieties}
\subsection{}
In this section, we will introduce another class of varieties associated with a quiver.
Since some examples of such varieties have appeared in the work of Dancer et al.,~\cite{Dancer:2013aa,Dancer:2013ta}, we call a variety of this class a DKS type variety.
The method to construct such varieties imitates the construction of the Nakajima quiver varieties.
The main difference is that we need to use the special linear group to replace the general linear group.

As we have said in \S~\ref{sub:comp-q}, one merit of the DKS type variety is that, with a constraint on the dimension vectors, the Weyl group of the quiver can act on a single DKS type variety, rather than just sending a variety to another one as Theorem~\ref{thm:nweyl}.
The rest of this section is devoted to establishing the existence of such a Weyl group action.

\subsection{DKS type varieties.}
\label{sub:dks}
Choose a quiver $Q$ as in Section~\ref{sec:quiver-varieties}.
About the vector spaces $V,W$, we now assume that $V_k,W_k$ are equipped with (complex) volume forms $\Pi_k,\Omega_k$ respectively.
In other words, we fix a basis vector $\Pi_k$ (resp.\ $\Omega_k$) for $\bigwedge^{\mathrm{top}}_{\bC} V^*_k$ (resp.\ $\bigwedge^{\mathrm{top}}_{\bC} W^*_k$), the determinant line of the dual of $V_k$ (resp.\ $W_k$).
Denote $(V_k,\Pi_k)_{k\in I}$ and $(W_k,\Omega_k)_{k\in I}$ to be $(V,\Pi)$ and $(W,\Omega)$ respectively.
Although we have introduced volume forms on vector spaces, we still define the space $\bfM(V,W)$ as in (\ref{eq:defm}).
The only difference is that from now on, an isomorphism between $V$ and $V'$ needs to preserve volume forms. And as before, we will fix a choice of $(W,\Omega)$ throughout the following construction.

The appearance of extra volume forms forces us to restrict the group $\rGL_{\bfv}$ to a subgroup $\rSL_{\bfv} \coloneqq \prod_{k\in I} \rSL(V_k,\Pi_k)$ of it.
Clearly, $\rSL_{\bfv}$ doesn't depend on the choice of volume forms.
In fact, it is the commutator subgroup of $\rGL_{\bfv}$.
$\rSL_{\bfv}$ also preserves symplectic form $\omega$ defined in (\ref{eq:sympformm}).
Denote the Lie algebra of $\rSL_{\bfv}$ by $\ksl_{\bfv}\coloneqq \oplus_{k\in I} \ksl(V_k)$.
Since the injection $\ksl_{\bfv} \hookrightarrow \kgl_{\bfv}$ induces a map $\kgl_{\bfv}^* \rightarrow \ksl_{\bfv}^*$, via the Killing form, we get a projection from $\kgl_{\bfv}$ to $\ksl_{\bfv}$.
The moment map $\mm$ of the $\rSL_{\bfv}$ action on $\bfM(V,W)$ is just the composition of the map $\mu$ defined in (\ref{eq:mm}) and such a projection.
More concretely, we have
\begin{equation}
  \label{eq:smm}
  \mm(B,i,j) = \big((\sum_{\sch} B_hB_{\bar{h}} - \sum_{\sct} B_{\bar{h}}B_{h} + i_kj_k)_0\big)_k \in \bigoplus_{k\in I} \ksl(V_k) = \ksl_{\bfv},
\end{equation}
where ${(\bullet)}_0$ means taking the trace-free part of a matrix.
Therefore, $\mm(B,i,j) = 0$ if and only if there exists $\lambda =\allowbreak (\lambda_1,\cdots,\lambda_n)\in \bC^n$ such that $\mu(B,i,j) = \lambda$, i.e.\ each component of $\mu(B,i,j)$ is a scalar matrix.\footnote{$\mm^{-1}(0)$ is isomorphic to the variety $\Lambda_{D,V}$ defined in~\cite[\S~3.1]{Lusztig:2000aa}.}

\begin{definition}
  \label{def:cM}
  Let $(V,\Pi),(W,\Omega)$ be vector spaces with volume forms attached to a quiver $Q$, whose dimension vectors are $\bfv,\bfw$.
  The affine GIT quotient defined using these data,
  \begin{equation}
    \label{eq:dks-v}
    \cM(\bfv,\bfw) \coloneqq \mm^{-1}(0) \dsl \rSL_{\bfv} = \Spec{\bC[\mm^{-1}(0)]^{\rSL_{\bfv}}},
  \end{equation}
  is called a \emph{DKS type variety}.
  Denote the quotient map from $\mm^{-1}(0)$ to $\cM(\bfv,\bfw)$ by $\pi$.
\end{definition}

When there is no ambiguity, we will write $\cM(\bfv,\bfw)$ as $\cM$ for short.

As before, we can also define $\cM$ as a \hkh\ quotient.
Let $\rSU_{\bfv}$ be the unitary subgroup of $\rSL_{\bfv}$ and $\mm_{\bR}$ be the real moment map for the $\rSU_{\bfv}$ action.
\begin{equation}
  \label{eq:dks-vr}
  \cM(\bfv,\bfw) = \big(\mm_{\bR}^{-1}(0) \cap \mm^{-1}(0) \big) / \rSU_{\bfv}.
\end{equation}

From the viewpoint of reduction in stages, $\cM(\bfv,\bfw)$ is just the intermediate product for the \hkh\ reduction of $\bfM(\bfv,\bfw)$ with respect to the $\rU_{\bfv}$ action.
More specifically, since $\cM$ inherits a \hkh\ Hamiltonian $T^n\coloneqq (S^1)^n$ action from $\bfM$, we can calculate the Nakajima quiver varieties with the parameter $\zeta = (\lambda,\theta)$ as follows,
\begin{equation*}
  \kM_{\zeta}(\bfv,\bfw) = \big(\ubar{\mu}^{-1}_{T^n}(\lambda) \cap \ubar{\mu}^{-1}_{T^n,\bR}(\theta)\big) / T^n,
\end{equation*}
where $\ubar{\mu}_{T^n}$ (resp.\ $\ubar{\mu}_{T^n,\bR}$) is the induced complex (resp.\ real) $T^n$ moment map on $\cM$.

As in \S~\ref{sub:cb}, $\cM$ can also be constructed from the quiver $Q^{\bfw}$.
For this purpose, we should fix a volume form on the vector space $V_{\infty}$ over the new vertex $\infty$.
Since $V_{\infty} \simeq \bC$, we can and will equip $V_{\infty}$ with the standard volume form.
Besides, in \S~\ref{sub:cb}, to compare the varieties constructed out of $Q$ and $Q^{\bfw}$, we have picked out a basis for each $W_i$.
For this new situation involving volume forms, we further require the chosen basis $\aset{x^i_1,\cdots,x^i_{w_i}}$ for $W_i$ satisfying $\Omega_i(x^i_1,\cdots,x^i_{w_i}) = 1$.
Now, the same arguments in \S~\ref{sub:cb} show that $\cM$ can be constructed by using either $Q$ or $Q^{\bfw}$.

\begin{remark*}
As we have mentioned, some examples of the DKS type variety have appeared in Dancer, Kirwan and Swann's work~\cite{Dancer:2013aa}.
  More precisely, in~\cite{Dancer:2013aa}, the authors have study the variety $\cM(\bfv,\bfw)$ associated with the quiver of $\rA_n$ type, where $\bfv = \tsp{(v_1,v_2,\cdots,v_n)}$, $\bfw = \tsp{(w_1,0,\cdots,0)}$ and $w_1 \ge v_1 \ge v_2 \ge \cdots \ge v_n \ge 0$.
  However, in the same work, the authors seem to prefer interpreting such special dimension vectors as partial flags of a vector space.
  Especially, they don't distinguish the roles played by $\bfw$ and $\bfv$.
  But in this paper, we are mostly interested in $\cM(\bfv,\bfw)$ satisfying $\bfw = \bfC\bfv$, which is not associated with a partial flag in general.
\end{remark*}

\subsection{The Weyl group action on $\cM$: the definition of $S_k$.}\label{sub:weyl-group-action}
In~\cite{Nakajima:1994aa}, to define the Weyl group action on $\aset{\kM_{\zeta}}$, Nakajima uses the gauge theoretic explanation of $\kM_{\zeta}$, c.f.~\cite{Kronheimer:1990aa}.
Later, several authors develop a more algebraic method to construct the same group action by using the so-called \emph{reflection functor},~\cite{Crawley-Boevey:1998aa,Lusztig:2000aa,Maffei:2002aa,Nakajima:2003aa}.
Here, we will try to modify the reflection functor method to fit our settings, which enables us to construct a $\bW$ action on $\cM$.

From now on, the dimension vectors are always assumed to satisfy $\bfw = \bfC \bfv$.
Moreover, we will work with the formulation of $\cM$ using $Q^{\bfw}$, which is a little more convenient.
The notation convention for the rest of this section is almost the same with that used in \S~\ref{sub:cb}.
Especially, it means that we use the same symbols for the corresponding moment maps on $\bfM(\hat{V})$ and $\bfM(V,W)$.
However, an exception is that unlike \S~\ref{sub:cb}, we will denote a point in $\bfM(\hat{V})$ by $B$ rather than by $\hat{B}$.
Since in the following we only consider points in $\bfM(\hat{V})$, we hope that such a simplified notation causes no ambiguity.

With these settings, we can restate Theorem~\ref{thm:main} in a more concrete way.
\begin{theorem}[{$=$Theorem~\ref{thm:main}}]
  \label{thm:main-re}
  The Weyl group $\bW$ of the quiver $Q$ of $\ade$ type can act on $\cM$ naturally. 
  Namely,
  \begin{enumerate}[label=(\arabic*)]
  \item for each simple reflection $s_k\in \bW$, $k\in I$, there is an automorphism $S_k$ on $\cM$ corresponding to $s_k$;
  \item $\aset{S_k}_{k\in I}$ satisfy the Coxeter relations (\ref{eq:cox1}--\ref{eq:cox3}).
  \end{enumerate}
\end{theorem}

As the statement of the above theorem hints, to define the Weyl group action on $\cM$, we proceed in two steps.
The \nth{1} step is to define an automorphism $S_k$ on $\cM$ for each simple reflection $s_k,k\in I$.
The \nth{2} step is to check that the defined $\aset{S_k}_{k\in I}$ satisfy the Coxeter relations.
Among these two steps, the \nth{2} step is rather lengthy and we will deal with it after \S~\ref{sub:def-sk}.
In the following three subsections, we concentrate on how to define $S_k$.
As a preparation for this task, we begin with investigating some properties of two special varieties, $\mm^{-1}(0)$ and $Z_k$.

\subsection{Properties of $\mm^{-1}(0)$.}
Recall that the variety $\mm^{-1}(0)$, as introduced in \S~\ref{sub:dks}, is the zero level set of the moment map $\mm$, (\ref{eq:smm}). 
To study properties of $\mm^{-1}(0)$, we use a result about the variety $\mu^{-1}(\lambda)$, which is a paraphrase of a theorem due to Crawley-Boevey in our situation.

Take $\alpha$ to be a dimension vector of $Q^{\bfw}$.
Let $p(\alpha) \coloneqq 1+ \sum_{h\in Q^{\bfw}}\alpha_{\tl{h}}\alpha_{\hd{h}} - \tsp{\alpha}\alpha$.
Especially, we have $p(\hat{\bfv}) = p\big((\bfv,1)\big) = \sum_{h\in Q}\bfv_{\tl{h}}\bfv_{\hd{h}} + \tsp{\bfv}\bfw - \tsp{\bfv}\bfv$.
For the definition of (positive) roots appearing in the following theorem, readers can find it in~\cite[\S~2]{Crawley-Boevey:2001aa}.
Here, we will not repeat this definition, since we do not use it in this paper.
\begin{theorem}[{\cite[Theorem~1.1 \& 1.2]{Crawley-Boevey:2001aa}}]
  \label{thm:cb}
For any $\lambda\in \bC^n$, let $\Lambda = \Lambda(\lambda) \coloneqq (\lambda,\allowbreak -\tsp{\bfv}\lambda)\in \bC^{n+1}$. The following statements are equivalent.
  \begin{enumerate}[label=(\arabic*)]
  \item\label{con:orb1} There is a closed orbit lying in $\mu^{-1}(\lambda)$ under the $\rGL_{\bfv}$ action and the isotropic subgroup of this orbit is trivial.
  \item\label{con:orb2}
$\hat{\bfv}$ is a positive root of $Q^{\bfw}$.
    Let $R_{\lambda} = \aset{\alpha\in \bZ_{\ge 0}^{I \sqcup \{\infty\}} | \alpha \text{ is a positive root of } Q^{\bfw} \allowbreak \text{and }\tsp{\Lambda} \alpha = 0}$.
    For any decomposition $\hat{\bfv} = \sum^r_{t =1} \beta^{(t)}$ with $r\ge 2$ and $\beta \in R_{\lambda}$, we have $p(\hat{\bfv}) > \sum^r_{t=1} p(\beta^{(t)})$.
  \end{enumerate}
  If one of the above conditions is satisfied, $\mu^{-1}(\lambda)$ is a reduced and irreducible complete intersection of dimension $\tsp{\bfv}\bfA\bfv + 2\tsp{\bfv}\bfw - \tsp{\bfv}\bfv$.
  Moreover, if $\lambda = 0$ satisfies above conditions, the moment map $\mu: \rM(\hat{\bfv}) \rightarrow \kgl_{\bfv}$ is a flat and surjective morphism.
\end{theorem}
As pointed by~\cite[p.~259]{Crawley-Boevey:2001aa}, by King's work on stability,~\cite{King:1994aa}, the existence of such a closed orbit is equivalent to an algebraic result, that is, the existence of a simple representation of the (deformed) preprojective algebra $\Pi_{0}$.
In fact, Crawley-Boevey chose to state his theorem in this algebraic form.
But as we don't use the preprojective algebra, we prefer the more geometric form of Theorem~\ref{thm:cb} as stated here.

Now we can show some crucial algebraic properties of $\mm^{-1}(0)$ as an application of the above theorem, which is the content of the following lemma.

\begin{lemma}
  \label{lm:mm}
  If $Q$ is of finite type, $\mm^{-1}(0)$ is a reduced and irreducible complete intersection of dimension $\tsp{\bfv}\bfA\bfv + 2\tsp{\bfv}\bfw - \tsp{\bfv}\bfv + n$.
\end{lemma}
\begin{proof}
  Recall our assumption on dimension vectors, i.e.\ $\bfw = \bfC \bfv$.
  As a result, we can use~\cite[Proposition~3.24 \& Corollary~10.8]{Nakajima:1998aa}, which guarantee the existence of an orbit in $\mu^{-1}(0)$ satisfying the condition~\ref{con:orb1} of Theorem~\ref{thm:cb} for a finite type quiver.
  Hence, Theorem~\ref{thm:cb} implies that $\mu:\rM(\hat{\bfv})\rightarrow \kgl_{\bfv}$ is a flat morphism.
  By using the base change with respect to the natural inclusion from $\bC^n$ to $\kgl_{\bfv}$, the restriction of $\mu: \mm^{-1}(0) \rightarrow \bC^n$ is also flat.
  On the other hand, we notice that for the set $R_{\lambda}$ appearing in the condition~\ref{con:orb2} of Theorem~\ref{thm:cb}, the following relation is true: $R_{\lambda} \subseteq R_0$ for any $\lambda \in \bC^n$.
  Since we have shown that $\lambda = 0$ satisfies the condition~\ref{con:orb2} of Theorem~\ref{thm:cb}, this relation about $R_{\lambda}$ implies the same condition is valid for any $\lambda \in \bC^n$.
  Therefore, $\mu^{-1}(\lambda)$ is reduced, irreducible and of dimension $\tsp{\bfv}\bfA\bfv + 2\tsp{\bfv}\bfw - \tsp{\bfv}\bfv$.
  Then, by the flatness of $\mu$ and~\cite[p.~137,~Proposition~3.8]{Liu:2002aa}, we conclude that $\mm^{-1}(0)$ is reduced and irreducible.
  Meanwhile, the dimension property of flat morphisms also implies the dimension of $\mm^{-1}(0)$ is $\tsp{\bfv}\bfA\bfv + 2\tsp{\bfv}\bfw - \tsp{\bfv}\bfv + n$, which in turn implies that $\mm^{-1}(0)$ is a complete intersection.
\end{proof}

\begin{remark*}
  In fact, with a fixed choice of $\bfv,\bfw$,~\cite[Corollary~10.8]{Nakajima:1998aa} enables us to relax the finiteness condition of $Q$ in the above lemma a little.
  Namely, let $L_{\bfw}$ be an integrable highest weight module of the Kac-Moody algebra associated with $Q$.
  Suppose that the highest weight of $L_{\bfw}$ is $\sum_{k=1}^nw_k\varpi_k$, where $\aset{\varpi_k}$ are fundamental weights.
  To ensure that Lemma~\ref{lm:mm} holds in this case, we should assume that $0$ occurs as a weight of $L_{\bfw}$.
  
  In particular, let the undirected graph of $Q$ be an affine Dynkin graph.
  For such a graph, there is a special vertex corresponding to the negative of the highest root of the underlying simple Lie algebra.
  Label this vertex with $0$.
  Choose $\bfv$ to be the unique vector such that $\bfC \bfv = 0$ and $v_0 = 1$.
  Then if we take $\bfw$ to be $0$, $Q$, $\bfv$, $\bfw$ will satisfy the condition in the last paragraph.
  This special non-finite example is particularly interesting because in this case $\cM(\bfv,\bfw)$ has a close relation with the ALE spaces, c.f.~\cite[p.372]{Nakajima:1994aa}.
\end{remark*}

From now on, we always assume that the quiver $Q$ appearing in this section is of finite type, i.e.\ $\ade$ type.

\subsection{The variety $Z_k$.}
We will define another variety $Z_k$ based on $\mm^{-1}(0)$, which is a variant of a variety introduced by Lusztig~\cite[\S~3.2]{Lusztig:2000aa}.
As in~\cite{Lusztig:2000aa}, $Z_k$ lies in the heart of the whole construction for the Weyl group action on $\cM$.
It is helpful to introduce some notations.

Fixing a vertex $k\in I$, we define the following space, c.f.~\cite[\S~3.2]{Lusztig:2000aa} and~\cite[Definition~27]{Maffei:2002aa},
\begin{equation*}
  T_k \coloneqq \bigoplus_{\hd{h} = k,h\in \DQw} V_{\tl{h}}.
\end{equation*}
For each $k\in I$, we fix a volume form for $T_k$ via choosing an order for elements in the set $\aset{V_{\tl{h}}|\hd{h} = k,h\in \DQw}$.
Moreover, for any $B\in \mm^{-1}(0)$, we can define two useful maps associated with $T_k$.
\begin{subequations}
  \begin{alignat}{2}
    \om_k(B) &\coloneqq \bigoplus_{{\hd{h} = k,h\in \DQ^{\bfw}}}B_{\bar{h}} &&:V_k \rightarrow T_k,\label{eq:def-om}\\
    \im_k(B) &\coloneqq \sum_{{\hd{h} = k,h\in \DQ^{\bfw}}} \epsilon(h)B_h &&:T_k \rightarrow V_k,\label{eq:def-im}\end{alignat}
\end{subequations}
where $\epsilon(h) = 1$ for $h\in Q^{\bfw}$ and $\epsilon(h) = -1$ for $h\in Q^{\bfw,\mathrm{op}}$.

Traditionally, e.g.~\cite{Lusztig:2000aa,Maffei:2002aa}, the maps $\om_k,\im_k$ defined as above are denoted by $a_k, b_k$ respectively.
However, with the $a_k, b_k$ notations, it is a little difficult to memorize the domain and the range of the maps denoted by $a_k,b_k$.
In view of this, we hope that our new notations can give readers some hints about the direction of the maps that they denote.
That is, $\om_k (B)$ (resp.\ $\im_k(B)$) denotes a map going \uline{out} of (resp.\ \uline{into}) $V_k$.

Let $\mu_k(B)$ be the $k^{\text{th}}$ component of $\mu(B)$.
Since $B\in \mm^{-1}(0)$, $\mu_k(B)$ is a scalar matrix.
As in \S~\ref{sub:quiver}, we will treat $\mu_k(B)$ as a scalar, that is, we will treat $\mu(B)$ as an element in $\bC^n$.
Bearing this in mind, by (\ref{eq:mm-v}), one can check that
\begin{equation}
  \label{eq:ba}
  \im_k(B)\om_k(B) = \mu_k(B)\id_{V_k}.
\end{equation}

Recall that our aim is to define an automorphism $S_k$ on $\cM$ corresponding a simple reflection $s_k$.
The following two open subvarieties of $\mm^{-1}(0)$ are quite relevant to this purpose,
\begin{align*}
  \Lambda_k^{\mathrm{s}} &\coloneqq \aset{B\in \mm^{-1}(0)|\,\text{$\im_k(B)$ is surjective}}, \\
  \Lambda_k^{\mathrm{i}} &\coloneqq \aset{B\in \mm^{-1}(0)|\,\text{$\om_k(B)$ is injective}},
\end{align*}
both of which are $\rSL_{\bfv}$ invariant.
Note that due to (\ref{eq:ba}), both $\Lambda_k^{\mathrm{s}}, \Lambda_k^{\mathrm{i}}$ contain an open subvariety $\aset{B\in \mm^{-1}(0) | \mu_k(B) \neq 0}$.
Besides, since $Q$ is of finite type, by using the definition of $\mu$, one can find points lying in $\Lambda_k^{\mathrm{s}}$ (resp.\ $\Lambda_k^{\mathrm{i}}$) but not in $\aset{B\in \mm^{-1}(0) | \mu_k(B) \neq 0}$.

\begin{definition}
  \label{def:zk}
  We define $Z_k$ to be a subvariety of $\mm^{-1}(0) \times \mm^{-1}(0)$ consisting of pairs $(B,B')$ satisfying the following conditions.
  \begin{enumerate}[label=C\arabic*.,ref=C\arabic*]
  \item\label{con:c1} $B_h = B'_h$ for any $h\in \DQw$ such that $\hd{h} \neq k$ and $\tl{h} \neq k$.
  \item\label{con:c2} The following short sequence is exact,
    \begin{equation}
      \label{eq:ses}
      \begin{tikzcd}
        0 \arrow[r] & V_k \arrow[r, "\om_k(B')"] & T_k  \arrow[r, "\im_k(B)"] & V_k  \arrow[r] & 0.
      \end{tikzcd}
    \end{equation}
    Moreover, the torsion of the above sequence is $1$ with respect to the volume forms of $V_k$ and $T_k$.
  \item\label{con:c3} $\om_k(B')\im_k(B') = \om_k(B)\im_k(B) - \mu_k(B)\id_{T_k}: T_k \rightarrow T_k$.
  \end{enumerate}
\end{definition}

For the torsion of a chain complex appearing in the condition~\ref{con:c2}, readers can find a good introduction on this topic in~\cite[Ch.~I]{Turaev:2001ab}.

\begin{remark}
  \label{rk:bb}
  We make some short comments on the conditions of $Z_k$.
  \begin{enumerate}[label=\Alph*.,ref=\Alph*]
  \item\label{rk:bb1} By~\cite[p.472, Remark]{Lusztig:2000aa} or~\cite[Lemma~28]{Maffei:2002aa}, for $(B,B')\in Z_k$, we have the equality $\mu(B') = s_k\big(\mu(B)\big)$, where $s_k$ acts on $\bC^n$ by (\ref{eq:weylac}).
  \item\label{rk:bb2} Furthermore, the results of Lusztig~\cite{Lusztig:2000aa} and Maffei~\cite{Maffei:2002aa} also show that the assumption that both $B,B'$ lie in $\mm^{-1}(0)$ can be weakened.
    In fact, if we only assume that $B\in \mm^{-1}(0)$ and $B'\in \bfM$, conditions~\ref{con:c1}-\ref{con:c3} imply that $B'$ falls into $\mm^{-1}(0)$.
  \item To ensure the existence of the short exact sequence (\ref{eq:ses}), we have used the assumption on the dimension vectors, i.e.\ $\bfw = \bfC \bfv$.
  \end{enumerate}
\end{remark}

By Definition~\ref{def:zk}, the following group acts naturally on $Z_k$.
\begin{equation*}
  G_{k,\bfv} = \rSL(V_k) \times \rSL(V_k) \times \prod_{l\neq k,l\in I} \rSL(V_l).
\end{equation*}
More precisely, for $g = (g_1,\cdots,g_k, g_k',\cdots,g_n)\in G_{k,\bfv}$, $(B,B')\in Z_k$ and $h\in \DQw$, the component of $g\cdot(B,B')$ at $h$ is
\begin{equation*}
  (g\cdot(B,B'))_h =
  \begin{cases}
    (g_{\hd{h}}B_hg_{\tl{h}}^{-1}, g_{\hd{h}}B'_hg_{\tl{h}}^{-1}), & \text{if } \hd{h}\neq k \text{ and } \tl{h}\neq k;\\
    (g_{\hd{h}}B_hg_k^{-1}, g_{\hd{h}}B'_h{g'_k}^{-1}), & \text{if } \tl{h} = k;\\
    (g_kB_hg_{\tl{h}}^{-1}, g'_kB'_hg_{\tl{h}}^{-1}), & \text{if } \hd{h} = k.
  \end{cases}
\end{equation*}
Via the injection of $\rSL(V_k)$ into the the first (resp.\ the second) copy of $\rSL(V_k)$ in the decomposition of $G_{k,\bfv}$, the above $G_{k,\bfv}$ action induces an $\rSL(V_k)$ action on $Z_k$, denoted by $\rho_1$ (resp.\ $\rho_2$).

Meanwhile, let $\opp$ (resp.\ $\opp'$) denote the projection from $\mm^{-1}(0) \times \mm^{-1}(0)$ to the first (resp.\ second) copy of $\mm^{-1}(0)$.
And the restriction of $\opp$ and $\opp'$ on $Z_k$ are denoted by the same symbols.
Similar to~\cite[Lemma~33]{Maffei:2002aa}, the following result holds.
\begin{lemma}
  \label{lm:pp}
  By the $\rSL(V_k)$ action $\rho_2$ (resp.\ $\rho_1$), the map $\opp:Z_k \rightarrow \Lambda_k^{\mathrm{s}}$ (resp.\ $\opp':Z_k \rightarrow \Lambda_k^{\mathrm{i}}$) yields a principal $\rSL(V_k)$ bundle structure on $Z_k$.
\end{lemma}
\begin{proof}
  Since two cases can be handled in a similar way, we only show one of them.
  We use the same notations as in Definition~\ref{def:zk}.
  Let $B\in \Lambda_k^{\mathrm{s}}$.
  We need to find $B'\in \mm^{-1}(0)$ such that $(B,B')\in Z_k$ and show that such a point $B'$ is unique up to an $\rSL(V_k)$ $\rho_2$ action.
  
  To show this result, we note the following fact: due to the torsion $1$ condition in the condition~\ref{con:c2}, for a fixed $\im_k(B)$, any two $\om_k(B')$ satisfying the condition~\ref{con:c2} coincide up to an $\rSL(V_k)$ $\rho_2$ action at most.
  Using this fact, we can reproduce the proof of~\cite[Lemma~33]{Maffei:2002aa} verbatim to show that for any $\om_k(B')$ satisfying the condition~\ref{con:c2}, there exists a unique $\im_k(B')$ such that the condition~\ref{con:c3} holds and the desired $B'$, as well as its uniqueness, is obtained consequently.
\end{proof}

\begin{remark*}
  We do not specify the topology used to define the principal bundle in the above lemma.
  This will not cause any ambiguities thanks to a result by Grothendieck,~\cite[Th\'eor\`eme~3]{Grothendieck:1958aa}, which asserts that when the structure group is a special linear group (as in our situation), the principal bundle is locally trivial with respect to any usual topology, e.g.\ Zariski topology, \'etale topology, etc.
  In fact, it is not difficult in our concrete case to find a local section on a Zariski open subset for $\opp$ or $\opp'$.
  Moreover, in view of Lemma~\ref{lm:mm}, the principal bundle structure implies that $Z_k$ is also reduced and irreducible.
\end{remark*}

\subsection{The definition of $S_k$.}
To define $S_k$, the variety $Z_k$ plays an intermediate role.
In fact, the following result is a direct consequence of Lemma~\ref{lm:pp}, which compares the invariant regular functions on several varieties.
\begin{proposition}
  \label{prop:rf}
  \begin{equation*}
    \label{eq:rf}
\bC[\Lambda_k^{\mathrm{s}}]^{\rSL_{\bfv}} \xrightarrow[\,\opp^{*}]{{\displaystyle \simeq}} \bC[Z_k]^{G_{k,\bfv}}  \xleftarrow[\opp'^*]{{\displaystyle \simeq}} \bC[\Lambda_k^{\mathrm{i}}]^{\rSL_{\bfv}}.
  \end{equation*}
\end{proposition}
Due to this proposition, $\opp'^{*,-1}\opp^{*}$ is an isomorphism from $\bC[\Lambda_k^{\mathrm{s}}]^{\rSL_{\bfv}}$ to $\bC[\Lambda_k^{\mathrm{i}}]^{\rSL_{\bfv}}$.
We denote this isomorphism by $\psi$ in the following.

\begin{propdef}
  \label{propdef:sk}
  For any simple reflection $s_k$, the following three varieties are naturally isomorphic to each other,
  \begin{equation}
    \label{eq:3-var}
    \Spec{\bC[\mm^{-1}(0)]^{\rSL_{\bfv}}} = \Spec{\bC[\Lambda_k^{\mathrm{s}}]^{\rSL_{\bfv}}} = \Spec{\bC[\Lambda_k^{\mathrm{i}}]^{\rSL_{\bfv}}}.
  \end{equation}
  Due to this fact, we define $S_k$ to be the automorphism of $\cM$ which makes the following diagram commute,
  \begin{equation*}
    \begin{tikzcd}
      \cM = &[-3em] \Spec{\bC[\mm^{-1}(0)]^{\rSL_{\bfv}}} \arrow[r, "S_k"] \arrow[d, equal]& \Spec{\bC[\mm^{-1}(0)]^{\rSL_{\bfv}}} \arrow[d, equal] &[-3em] = \cM\\
      & \Spec{\bC[\Lambda_k^{\mathrm{s}}]^{\rSL_{\bfv}}} \arrow[r, "(\Spec{\psi})^{-1}"]& \Spec{\bC[\Lambda_k^{\mathrm{i}}]^{\rSL_{\bfv}}} &
    \end{tikzcd}.
  \end{equation*}
\end{propdef}
\begin{proof}
  By Proposition~\ref{prop:rf}, we know that
  \begin{equation*}
    \Spec{\psi}: \Spec{\bC[\Lambda_k^{\mathrm{i}}]^{\rSL_{\bfv}}} \longrightarrow \Spec{\bC[\Lambda_k^{\mathrm{s}}]^{\rSL_{\bfv}}}
  \end{equation*}
  is an isomorphism.
  Hence, to show that $S_k$ in the above definition is an automorphism on $\cM$, we only need to show (\ref{eq:3-var}).
  
  Since $\mm^{-1}(0)$ is a complete intersection by Lemma~\ref{lm:mm}, $\mm^{-1}(0)$ is Cohen-Macaulay, especially satisfying Serre's S$_2$ property.
  As we have noted, $\Lambda_k^{\mathrm{s}}$ (resp.\ $\Lambda_k^{\mathrm{i}}$) contains but not equals to $\aset{B\in \mm^{-1}(0) | \mu_k(B) \neq 0}$.
  In other words, $\mm^{-1}(0) - \Lambda_k^{\mathrm{s}}$ (resp.\ $\mm^{-1}(0) - \Lambda_k^{\mathrm{i}}$) is a subvariety strictly contained in $\mm^{-1}(0) \cap \mu_k^{-1}(0)$.
  Therefore, $\mm^{-1}(0) - \Lambda_k^{\mathrm{s}}$ (resp.\ $\mm^{-1}(0) - \Lambda_k^{\mathrm{i}}$) is of codimension at least $2$ in $\mm^{-1}(0)$.
  As a result, we can use~\cite[Th\'eor\`eme~(5.10.5)]{Grothendieck:1965aa}, which asserts that a regular function can be extended from $\Lambda_k^{\mathrm{s}}$ (resp.\ $\Lambda_k^{\mathrm{i}}$) to $\mm^{-1}(0)$, that is, $\bC[\mm^{-1}(0)] = \bC[\Lambda_k^{\mathrm{s}}] = \bC[\Lambda_k^{\mathrm{i}}]$.
Consequently, we also have $\bC[\mm^{-1}(0)]^{\rSL_{\bfv}} = \bC[\Lambda_k^{\mathrm{s}}]^{\rSL_{\bfv}} = \bC[\Lambda_k^{\mathrm{i}}]^{\rSL_{\bfv}}$, from which (\ref{eq:3-var}) follows.
\end{proof}

\begin{remark}
  \label{rk:sk}
  Taking $B\in \Lambda_k^{\mathrm{s}}$ and $(B,B')\in Z_k$, the above construction implies $S_k(\pi(B)) = \pi(B')$.
Especially, by Remark~\ref{rk:bb}, we know that $S_k$ preserves the open and dense subvariety $\aset{\pi(B)| \mu_k(B) \ne 0} \subseteq \cM$.
\end{remark}

\subsection{The Weyl group action on $\cM$: the Coxeter relations.}\label{sub:def-sk}

In this subsection, we finish the proof of Theorem~\ref{thm:main-re} by showing that $\aset{S_k}_{k\in I}$ defined in previous subsection indeed satisfy the Coxeter relation~(\ref{eq:cox1}-\ref{eq:cox3}).

Technically, for fixed $k,l\in I$, we use the following open and dense subvariety of $\cM$,
\begin{equation*}
  \cM_{kl} \coloneqq \aset{\pi(B)| \mu_k(B)\neq 0, \mu_l(B)\neq 0, \mu_k(B)+\mu_l(B)\neq 0},
\end{equation*}
whose preimage in $\mm^{-1}(0)$ is denoted by $\bfN_{kl}$.
One can check that $S_k$ and $S_l$ preserve $\cM_{kl}$.
As a result, to show that $S_k$ and $S_l$ satisfy the Coxeter relation, it is sufficient to verify that the same is true for $S_k|_{\cM_{kl}}$ and $S_l|_{\cM_{kl}}$.
Therefore, in the rest of this section, for any $B\in \mm^{-1}(0)$, we can and will further assume that $B\in \bfN_{kl}$, which causes no loss of generality.
For simplicity, we denote $S_k|_{\cM_{kl}}$ and $S_l|_{\cM_{kl}}$ by $S_k$ and $S_l$ for short.

In the following discussion, given $(B,B')\in Z_k$, symbols $\om_k(B),\im_k(B)$ (resp.\ $\om_k(B'),\im_k(B')$) appear frequently.
Thus, we will write them as $\om_k,\im_k$ (resp.\ $\om'_k,\im'_k$) for short.

Moreover, as our vector spaces are always equipped with volume forms, we often need to take a basis of a vector space compatible with its volume form.
For a vector space with a volume form $(Y,\Omega)$, a basis $\aset{y_i}$ of $Y$ (or a multivector in top degree $\wedge_{i=1}^{\dim{Y}}y_i$) is called \emph{normed} if $\Omega(y_1,\cdots,y_{\dim{V}}) = 1$.

Of the three Coxeter relations, (\ref{eq:cox1}) and (\ref{eq:cox2}) are relatively straightforward to verify and (\ref{eq:cox3}) needs more efforts.
Let us begin with the easy part.

\subsection{Coxeter relation (\ref{eq:cox1}) \& (\ref{eq:cox2}).} These two relations can be verified by direct calculation.
\paragraph{Check {(\ref{eq:cox1})}.} To check (\ref{eq:cox1}), we are going to show that if $(B,B')\in Z_k$, then the same is true for $(B',B)\in Z_k$.

Among the three conditions in Definition~\ref{def:zk}, the condition~\ref{con:c1} and~\ref{con:c3} hold trivially for $(B',B)$.
And by the result of~\cite[Remark~36]{Maffei:2002aa}, we know that the chain complex (\ref{eq:ses}) in the condition~\ref{con:c2} is also exact for $(B',B)$.
All we need to verify is that the torsion of (\ref{eq:ses}) is $1$ when $\om_k',\im_k$ are substituted by $\om_k,\im_k'$ respectively.

Now, choose a normed basis $\aset{x^k_{\alpha}}_{\alpha=1}^{v_k}$ of $V_k$.
For $B \in \bfN_{kl}$, by (\ref{eq:ba}), one has
\begin{equation}
  \label{eq:2.1a1}
  \im_k\big(\om_k(x^k_{\alpha}/\lambda_k)\big) = x^k_{\alpha}
\end{equation}
where $\alpha= 1,\cdots,v_k$ and $\lambda_k = \mu_k(B)$.
On the other hand, by Remark~\ref{rk:bb}~\ref{rk:bb1}, we can show
\begin{equation}
  \label{eq:2.1a2}
  \im'_k\big(\om'_k(-x^k_{\alpha}/\lambda_k)\big) = x^k_{\alpha}
\end{equation}
in a similar way.

As the torsion of (\ref{eq:ses}) for $(B,B')$ is $1$, by (\ref{eq:2.1a1}) and the definition of torsion,~\cite[Definition~1.2]{Turaev:2001ab}, we know that $\big(\wedge_{\alpha=1}^{v_k} \om'_k(x^k_{\alpha})\big) \wedge \big(\wedge_{\alpha=1}^{v_k} \om_k(x^k_{\alpha}/\lambda_k)\big)$ is normed for $T_k$.
Since $v_k^2 + v_k$ is always even, we have
\begin{equation*}
  \big(\wedge_{\alpha=1}^{v_k} \om'_k(x^k_{\alpha})\big) \wedge \big(\wedge_{\alpha=1}^{v_k} \om_k(x^k_{\alpha}/\lambda_k)\big) = \big(\wedge_{\alpha=1}^{v_k} \om_k(x^k_{\alpha})\big) \wedge \big(\wedge_{\alpha=1}^{v_k} \om'_k(-x^k_{\alpha}/\lambda_k)\big).
\end{equation*}
Hence, the r.h.s.\ of above equality is also normed for $T_k$, which, by (\ref{eq:2.1a2}), implies the torsion of (\ref{eq:ses}) is also equal to $1$ for $(B',B)$.

\paragraph{Check {(\ref{eq:cox2})}.} This relation is a direct consequence of the definition of $Z_k$.

To see this, we take $(B,B')\in Z_k$, $(B'',B')\in Z_l$ and $(B^{(1)},B)\in Z_l$.
Since $\bfA_{kl}= 0$, by the condition~\ref{con:c1} in Definition~\ref{def:zk}, we have $(B'',B^{(1)})\in Z_k$.
Hence, by Remark~\ref{rk:sk}, one can see that $S_kS_l(\pi(B)) = \pi(B'') = S_lS_k(\pi(B))$ and the Coxeter relation (\ref{eq:cox2}) follows.

\subsection{Some linear algebra results.}
Before checking (\ref{eq:cox3}) for $S_k$, we need some linear algebra results as a preparation.
In essence, the results given here are minor extension of the corresponding results in~\cite{Maffei:2002aa} by taking the torsion into consideration.

\begin{lemma}
  \label{lm:la1}
  Let $V,W,X,Y,Z$ be finite-dimensional vector spaces with volume forms and $\alpha,\beta,\gamma,\delta,\epsilon,\varphi$ be linear maps between them as indicated in the diagram below.
  The chain complex,
  \begin{equation}
    \label{eq:la1cx1}
    \begin{tikzcd}[ampersand replacement=\&]
      0 \arrow[r] \& V \arrow[r, "L_1 = {\renewcommand\arraystretch{0.7}\begin{pmatrix}
        \alpha \\ \beta \\ \gamma
      \end{pmatrix}}"] \&[1.2em] W \oplus X \oplus Y \arrow[r, "L_2 = {\renewcommand\arraystretch{0.7} \setlength\arraycolsep{2pt}\begin{pmatrix}
  \delta & 0 & -1 \\ 0 & \epsilon & \varphi
\end{pmatrix}}"] \&[2.2em] Y \oplus Z \arrow[r] \& 0,
    \end{tikzcd}
  \end{equation}
  is exact if and only if the chain complex,
  \begin{equation}
    \label{eq:la1cx2}
    \begin{tikzcd}[ampersand replacement=\&]
      0 \arrow[r] \& V \arrow[r, "{\renewcommand\arraystretch{0.7}\begin{pmatrix}
          \alpha \\ \beta
        \end{pmatrix}}"] \& W \oplus X \arrow[r, "{\renewcommand\arraystretch{0.7}\setlength\arraycolsep{2pt}\begin{pmatrix}
        \varphi\delta & \epsilon
      \end{pmatrix}}"] \&[1em] Z \arrow[r] \& 0,
    \end{tikzcd}
  \end{equation}
  is exact and $\gamma = \delta\alpha$. And the torsions of two chain complexes coincide up to a sign $(-1)^{\dim{Z}\dim{Y} + \dim{Y}}$.
\end{lemma}
\begin{proof}
  We only show how to calculate the torsion, which is the only difference between the above lemma and~\cite[Lemma~40]{Maffei:2002aa}.
  It is convenient to assume that the torsion of the complex (\ref{eq:la1cx2}) in the above lemma is $1$.
  This assumption causes no loss of generality because the general case can always be reduced to this case by rescaling the volume form on $V$.
  
  Choose normed top degree multivectors $\wedge_iv_i$, $\wedge_jz_j$ of $V$ and $Z$ respectively\footnote{Unfortunately, the symbol $v_i$ has been used for a component of a dimension vector.
    The good news is that we only use $v_i$ to denote a vector in $V$ in this subsection.
    Outside this subsection, $v_i$ always denotes a component of a dimension vector.}.
  And let $\wedge_j(w_j+ x_j)$ be a preimage of $\wedge_jz_j$ in $W \oplus X$.
  Since we have assumed the torsion of the complex (\ref{eq:la1cx2}) is $1$, $\big(\wedge_i(\alpha v_i+ \beta v_i)\big) \wedge \big(\wedge_j(w_j+ x_j)\big)$ is normed for $W \oplus X$.
  
  Now we also choose a normed multivector $\wedge_ky_k$ for $Y$.
  Then, the following equality holds because $\aset{z_k}$ is a basis for $Z$,
  \begin{equation}
    \label{eq:zy}
    (\wedge_jz_j) \wedge \big(\wedge_k (-y_k+ \varphi y_k)\big) = (-1)^{\dim{Z}\dim{Y} + \dim{Y}} (\wedge_ky_k) \wedge (\wedge_jz_j).
  \end{equation}
  We can choose the preimages of $z_j$ and $-y_k+ \varphi y_k$ in $W \oplus X \oplus Y$ as follows,
  \begin{equation}
    \label{eq:wz}
    L_2(w_j,x_j,\delta w_j) = \epsilon x_j + \varphi\delta w_j = z_j,\quad L_2(0,0,y_k) = -y_k+ \varphi y_k.
  \end{equation}
  Since $\aset{y_k}$ spans $Y$, one also has
  \begin{multline}
    \label{eq:wzy}
    \big(\wedge_i(L_1(v_i))\big) \wedge \big(\wedge_j(w_j + x_j + \delta w_j) \big) \wedge (\wedge_k y_k) \\
    = \big(\wedge_i(\alpha v_i + \beta v_i)\big) \wedge \big(\wedge_j(w_j + x_j) \big) \wedge (\wedge_k y_k),
  \end{multline}
  the r.h.s.\ of which is normed for $W \oplus X \oplus Y$ by our choice of $\aset{v_i}$, $\aset{w_j+x_j}$ and $\aset{y_k}$.
  
  Hence, by (\ref{eq:zy}), (\ref{eq:wz}), (\ref{eq:wzy}), we know that the torsion of the complex (\ref{eq:la1cx1}) in the lemma is $(-1)^{\dim{Z}\dim{Y} + \dim{Y}}$.
  The proof of the lemma is completed.
\end{proof}

We also need another linear algebra lemma of the same nature.
\begin{lemma}
  \label{lm:la2}
  Let $U,V,W,X,Y,Z$ be finite-dimensional vector spaces with volume forms and $\alpha,\beta,\gamma,\delta,\phi,\rho,\eta$ be linear maps between them as indicated in the diagrams below.
  Suppose that the following two chain complexes are exact and of torsion $\tau_1$, $\tau_2$ respectively,
  \begin{subequations}
    \begin{equation}
      \begin{tikzcd}[ampersand replacement=\&]
        0 \arrow[r] \& V \arrow[r, "{\renewcommand\arraystretch{0.7}\begin{pmatrix}
            \alpha \\ \beta
          \end{pmatrix}}"] \&[1.2em] W \oplus U \arrow[r, "{\renewcommand\arraystretch{0.7} \setlength\arraycolsep{2pt}\begin{pmatrix}
            \gamma & -\delta 
          \end{pmatrix}}"] \&[1.2em] Y  \arrow[r] \& 0,   
      \end{tikzcd}\label{eq:cx1}
    \end{equation}
    \begin{equation}
      \begin{tikzcd}[ampersand replacement=\&]
        0 \arrow[r] \& U \arrow[r, "{\renewcommand\arraystretch{0.7}\begin{pmatrix}
            \phi \\ \delta
          \end{pmatrix}}"] \&[1.2em] X \oplus Y \arrow[r, "{\renewcommand\arraystretch{0.7} \setlength\arraycolsep{2pt}\begin{pmatrix}
            \rho & \eta 
          \end{pmatrix}}"] \&[1.2em] Z \arrow[r] \& 0.
      \end{tikzcd}\label{eq:cx2}
    \end{equation}
  \end{subequations}
  Then, the diagram
  \begin{equation}
    \label{eq:cx3}
    \begin{tikzcd}[ampersand replacement=\&]
      0 \arrow[r] \& V \arrow[r, "{\renewcommand\arraystretch{0.7}\begin{pmatrix}
          \alpha \\ \phi\beta \\ \delta\beta
        \end{pmatrix}}"] \&[1.2em] W \oplus X \oplus Y \arrow[r, "L = {\renewcommand\arraystretch{0.7} \setlength\arraycolsep{2pt}\begin{pmatrix}
          \gamma & 0 & -1 \\ 0 & \rho & \eta
        \end{pmatrix}}"] \&[2.2em] Y \oplus Z  \arrow[r] \& 0
    \end{tikzcd}
  \end{equation}
  is an exact short sequence, the torsion of which is equal to $\tau_1\tau_2$.
\end{lemma}
\begin{proof}
  As before, we only need to show the lemma with the assumption that torsions of (\ref{eq:cx1}), (\ref{eq:cx2}) are equal to $1$.
  We choose $\wedge_iv_i$, $\wedge_jy_j$, $\wedge_kz_k$ to be normed top degree multivectors of $V$, $Y$ and $Z$ respectively.

  We first show the exactness of (\ref{eq:cx3}).
  Since (\ref{eq:cx1}) and (\ref{eq:cx2}) are chain complexes, we can see that diagram (\ref{eq:cx3}) is a chain complex.
  Now, we are going to show $L$ is surjective, that is, (\ref{eq:cx3}) is exact at place $Y \oplus Z$.
  Due to the definition of $L$, we only need to find the preimage of $z_k$ in $W \oplus X \oplus Y$.

  By the exactness of (\ref{eq:cx2}), we know that there exist $\bar{x}_k\in X$, $\bar{y}_k\in Y$ such that $\rho\bar{x}_k + \eta\bar{y}_k = z_k$.
  And the exactness of (\ref{eq:cx1}) enables us to find $\dot{w}_k\in W$, $\bar{u}_k\in U$ such that $\gamma\dot{w}_k - \delta\bar{u}_k = \bar{y}_k$.
  Let $\dot{x}_k = \bar{x}_k + \phi\bar{u}_k$, $\dot{y}_k = \bar{y}_k + \delta\bar{u}_k$.
  One can check that
  \begin{equation}
    \label{eq:z-pi}
    L(\dot{w}_k,\dot{x}_k,\dot{y}_k) = z_k.
  \end{equation}
  
Using a similar argument, we can show the exactness (\ref{eq:cx3}) at other two places.
  The details are left to readers.

  The next step is calculating the torsion of (\ref{eq:cx3}).
  By using (\ref{eq:cx1}) again, we can find $w_j'\in W$, $u_j'\in U$ such that $\gamma w_j' - \delta u_j' = y_j$, which implies that
  \begin{equation}
    \label{eq:y-pi}
    L(w_j',\phi u_j', \delta u_j') = y_j + (\rho\phi u_j' + \eta\delta u_j').
  \end{equation}
  Since $(\rho\phi u_j' + \eta\delta u_j') \in Z$, by (\ref{eq:z-pi}) and (\ref{eq:y-pi}), we know that
  \begin{multline}
    \label{eq:yz}
    \big(\wedge_jL(w_j',\phi u_j', \delta u_j')\big) \wedge \big(\wedge_kL(\dot{w}_k,\dot{x}_k,\dot{y}_k)\big) \\
    = \big(\wedge_j(y_j + \rho\phi u_j' + \eta\delta u_j')\big) \wedge (\wedge_k z_k) = (\wedge_j y_j) \wedge (\wedge_k z_k)
  \end{multline}
  is normed in $Y \oplus Z$.
  Moreover, recall that the torsion of (\ref{eq:cx1}) is $1$, which means that
  \begin{equation}
    \label{eq:wu}
    \big(\wedge_i(\alpha v_i + \beta v_i)\big) \wedge \big(\wedge_j(w_j' + u_j')\big)
  \end{equation}
  is normed for $W \oplus U$.

  On the other hand, the exactness and the torsion $1$ property of (\ref{eq:cx2}) imply that the following chain complex
  \begin{equation*}
    \begin{tikzcd}[ampersand replacement=\&]
      0 \arrow[r] \& W \oplus U \arrow[r, "{\renewcommand\arraystretch{0.7}\begin{pmatrix}
          1 \\ \phi \\ \delta
        \end{pmatrix}}"] \&[1.2em] W \oplus X \oplus Y \arrow[r, "{\renewcommand\arraystretch{0.7} \setlength\arraycolsep{2pt}\begin{pmatrix}
          0 & \rho & \eta 
        \end{pmatrix}}"] \&[1.2em] Z \arrow[r] \& 0
    \end{tikzcd}
  \end{equation*}
  is exact and of torsion $1$.
  Since the multivector (\ref{eq:wu}) is normed for $W \oplus U$, by using (\ref{eq:z-pi}), the torsion $1$ property of the above chain complex implies that the multivector
  \begin{equation}
    \label{eq:wxy}
    \big(\wedge_i(\alpha v_i + \phi\beta v_i + \delta\beta v_i)\big) \wedge \big(\wedge_j(w_j' + \phi u_j' + \delta u_j')\big) \wedge \big(\wedge_k(\dot{w}_k + \dot{x}_k + \dot{y}_k)\big)
  \end{equation}
  is normed for $W \oplus X \oplus Y$.
  Now, since both (\ref{eq:yz}) and (\ref{eq:wxy}) are normed multivectors, we can conclude that the torsion of (\ref{eq:cx3}) is equal to $1$.
\end{proof}

The following corollary is an immediate consequence of the above two lemmas. 
\begin{corollary}
  \label{cor:la}
  With the same assumptions and notations as in Lemma~\ref{lm:la2}, the following chain complex is exact, 
  \begin{equation}
    \label{eq:cla}
    \begin{tikzcd}[ampersand replacement=\&]
      0 \arrow[r] \& V \arrow[r, "{\renewcommand\arraystretch{0.7}\begin{pmatrix}
          \alpha \\ \phi\beta
        \end{pmatrix}}"] \& W \oplus X \arrow[r, "{\renewcommand\arraystretch{0.7}\setlength\arraycolsep{2pt}\begin{pmatrix}
        \eta\gamma & \rho
      \end{pmatrix}}"] \&[1em] Z \arrow[r] \& 0,
    \end{tikzcd}
  \end{equation}
  and the torsion of this chain complex is $(-1)^{\dim{Z}\dim{Y} + \dim{Y}}\tau_1\tau_2$.
\end{corollary}
\begin{proof}
  Since (\ref{eq:cx1}) is a chain complex, we have $\gamma\alpha = \delta\beta$.
  Then we can use Lemma (\ref{lm:la1}) to calculate the torsion of (\ref{eq:cla}) from the torsion of (\ref{eq:cx3}).
\end{proof}
\begin{remark*}
  It is easy to see that if (\ref{eq:cx2}) and (\ref{eq:cx3}) are exact, then (\ref{eq:cx1}) is a chain complex and exact.
  Moreover, if we already know that (\ref{eq:cx1}) is a chain complex, then the exactness of (\ref{eq:cx2}) and (\ref{eq:cla}) can also yield the exactness of (\ref{eq:cx1}).
\end{remark*}

\subsection{Coxeter relation (\ref{eq:cox3}).}
With the above preparation, we can verify $\aset{S_k}$ satisfying the Coxeter relation (\ref{eq:cox3}).
From now on, we assume that $\bfA_{kl} = 1$.

We introduce two auxiliary subspaces of $\bfN_{kl}\times \bfN_{kl}$,
\begin{gather*}
  \begin{multlined}
    Z_{klk} \coloneqq \aset{(B,B^{\pt}) \in \bfN_{kl}\times \bfN_{kl} | \exists\, B',B''\in \bfN_{kl}\; \text{such that}\, (B,B')\in Z_k,\\
      (B',B'')\in Z_l\;\text{and}\;(B'',B^{\pt})\in Z_k },
  \end{multlined}\\
  \begin{multlined}
    Z_{lkl} \coloneqq \aset{(B,B^{\pt}) \in \bfN_{kl}\times \bfN_{kl} | \exists\, B',B''\in \bfN_{kl}\; \text{such that}\, (B,B')\in Z_l,\\
      (B',B'')\in Z_k\;\text{and}\;(B'',B^{\pt})\in Z_l }.
  \end{multlined}
\end{gather*}
Note that $S_kS_lS_k(\pi(B)) = \pi(B^{\pt})$ (resp. $S_lS_kS_l(\pi(B)) = \pi(B^{\pt})$) if and only if $(B,B^{\pt}) \in Z_{klk}$ (resp. $(B,B^{\pt}) \in Z_{lkl}$).
Therefore, to show (\ref{eq:cox3}), it is necessary and sufficient to show the following result.
\begin{proposition}
  \label{prop:klk}
  If $\bfA_{kl} = 1$, then we have $Z_{klk} = Z_{lkl}$.
\end{proposition}

As before, if we forget the torsion condition in Definition~\ref{def:zk}, the method to show the above proposition has appeared in~\cite[p.680]{Maffei:2002aa}.
Hence, we need to check the argument carefully to ensure that the same method still works even with the extra constraint on the torsion.

\begin{proof}
  We only need to show
  \begin{equation}
    \label{eq:one-in}
    Z_{klk}\subseteq Z_{lkl}
  \end{equation}
  because by exchanging positions of $k$ and $l$ we can obtain the inclusion in the opposite direction in the same way.

  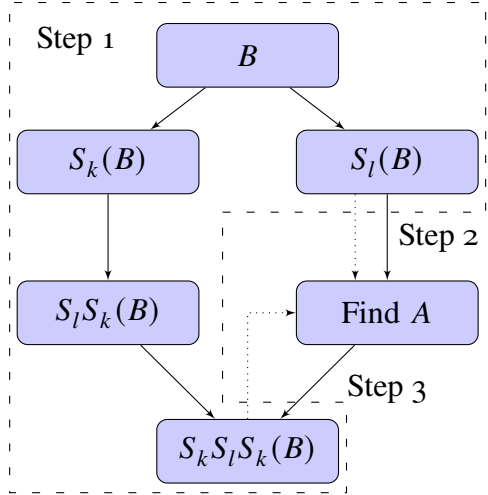
\begin{figure}[htbp]
    \centering
    \begin{tikzpicture}[node distance = 2cm, auto]
      \node [block] (init) {$B$};
      \node [block, below left of=init, xshift=-1em] (l1) {$S_k(B)$};
      \node [block, below right of=init, xshift=1em] (r1) {$S_l(B)$};
      \node [block, below of=l1] (l2) {$S_lS_k(B)$};
      \node [block, below of=r1] (r2) {Find $A$};
      \node [block, below right of=l2, xshift=1em, yshift=-1em] (fin) {$S_kS_lS_k(B)$};

      \node (text) at (-2.25,0.2) {Step 1};

      \path [line] (init) --  (l1);
      \path [line] (init) --  (r1);
      \path [line] (l1) --  (l2);
      \path [line] (r1) -- node {Step 2}(r2);
      \path [line,dotted,transform canvas={xshift=-1em}] (r1) -- (r2);
      \path [line] (l2) --  (fin);
      \path [line] (r2) -- node [near start] {Step 3} (fin);
      \path [line,dotted] (fin) |- (r2);

      \draw [loosely dashed](-0.33,-2.08) -- (3.13,-2.08) -- (3.13,0.7) -- (-3.13,0.7) -- (-3.13, -5.8) -- (1.3,-5.8) -- (1.3,-4.6) -- (-0.33,-4.6) -- cycle;
    \end{tikzpicture}

\caption{Idea of the proof of (\ref{eq:one-in}).}\label{fig:idea}
  \end{figure}

  By the definition of $Z_{klk}$, for $(B,B^{\pt})\in Z_{klk}$, we can find $B',B''\in \bfN_{kl}$ such that $(B,B')\in Z_k$, $(B',B'')\in Z_l$ and $(B'',B^{\pt})\in Z_k$.
However, instead of using symbols $B'$ (resp.\ $B''$, $B^{\pt}$), we will use more self-explanatory symbols $S_k(B)$ (resp.\ $S_lS_k(B)$, $S_kS_lS_k(B)$) to denote this intermediate object.

  The idea of the proof can be explained using Figure~\ref{fig:idea}.
  The \nth{1} step is to write down the objects encompassed by the dashed line, i.e.\ $B$, $S_k(B)$, \ldots, $S_kS_lS_k(B)$, explicitly and find the relations between them using Definition~\ref{def:zk}.
  The \nth{2} step is to use the information of $S_l(B)$ and $S_kS_lS_k(B)$ to find a point $A\in \mm^{-1}(0)$ such that $(S_l(B),A)\in Z_k$.
  The last step is to check that $(A, S_kS_lS_k(B))\in Z_l$.

  \paragraph{Step 1.} Let us start with the analysis of $B$, $S_k(B)$, \ldots, $S_kS_lS_k(B)$.

  Let $\mu_k(B) = \lambda_k$, $k\in I$.
  Since $\bfA_{kl} = 1$, there exists a unique edge $h_0\in \DQw$ going from $l$ to $k$.
  Set
  \begin{equation*}
    \varepsilon \coloneqq \epsilon(h_0).
  \end{equation*}

  In the following, we use two subspaces of $T_k$ and $T_l$ frequently.
  \begin{equation*}
    R_k \coloneqq \oplus_{h\in \DQw,\hd{h}= k,\tl{h}\neq l} V_{\tl{h}},\quad R_l \coloneqq \oplus_{h\in \DQw,\hd{h}= l,\tl{h}\neq k} V_{\tl{h}}.
  \end{equation*}
Recall that the volume form on $T_k$ is determined by a permutation of the set $\aset{V_{\tl{h}}|h\in \DQw,\hd{h}= k}$.
  Such a permutation also induces a volume form on $R_k$. (For a permutation of the set $\aset{V_{\tl{h}}|h\in \DQw,\hd{h}= k}$, we delete the element $V_l$.)
  With this volume form on $R_k$, the canonical isomorphism between $R_k \oplus V_l$ and $T_k$ preserves the volume form.
  We choose the volume form for $R_l$ in a similar way.

  Due to Definition~\ref{def:zk} and $\bfA_{kl}=1$, when comparing two objects $B_1$ and $B_2$ in neighboring positions in Figure~\ref{fig:idea}, for example, $B_1 = B$ and $B_2 = S_k(B)$, the only difference between $B_1$ and $B_2$ comes from the components corresponding to maps of following possibilities: between $R_k$ and $V_k$, between $V_k$ and $V_l$, or between $V_l$ and $R_l$.
Therefore, we will represent each object in Figure~\ref{fig:idea} by a diagram which highlights the ``key'' components of this object.

  For example, we represent $B$ by the following diagram,
  \begin{equation}
    \label{eq:ab1}
    B:\;
    \begin{tikzcd}
      R_k \arrow[r, rightharpoonup, shift left = 0.3ex, "d_k"] & \arrow[l, rightharpoonup, shift left = 0.3ex, "c_k"] V_k[\lambda_k] \arrow[r, rightharpoonup, shift left = 0.3ex, "\alpha"] & \arrow[l, rightharpoonup, shift left = 0.3ex, "\beta"] V_l[\lambda_l] \arrow[r, rightharpoonup, shift left = 0.3ex, "c_l"] & \arrow[l, rightharpoonup, shift left = 0.3ex, "d_l"] R_l.
    \end{tikzcd}
  \end{equation}
  Here and in the similar diagram appearing later, we use the notation convention as follows.
  \begin{enumerate}[label=(\arabic*)]
  \item The symbol $B$ at leftmost indicates that this diagram represents the object $B$ in Figure~\ref{fig:idea}.
  \item The symbol $V_k[\lambda_k]$ (resp.\ $V_l[\lambda_l]$) just denotes the vector space $V_k$ (resp.\ $V_l$), where the additional bracket is only used to mark the value of $\mu_k(B)$ (resp.\ $\mu_l(B)$).
  \item The symbol $\alpha$ (resp.\ $\beta$) in diagram (\ref{eq:ab1}) is defined to be the map $B_{\bar{h}_0}$ (resp.\ $B_{h_0}$).
  \item Morevover, the symbol $c_k$ (resp.\ $d_k$) denotes the map which is the composition of $\om_k(B)$ (resp.\ $\im_k(B)$) and the projection from $T_k$ to $R_k$ (resp.\ the injection from $R_k$ to $T_k$).
    And $c_l$, $d_l$ are defined in a similar way by replacing $k$ with $l$.
  \end{enumerate}

  With this notation convention, we will write down the diagrams for $S_k(B)$, $S_lS_k(B)$, $S_kS_lS_k(B)$ and $S_l(B)$ one by one.

\subparagraph{The diagram for ${S_k(B)}$.}
  Since $(B,S_k(B))\in Z_k$, by Remark~\ref{rk:bb}~\ref{rk:bb1}, we have $\mu_k(S_k(B)) = -\lambda_k$ and $\mu_l(S_l(B)) = \lambda_l + \lambda_k$.
  Hence, we have the diagram
  \begin{equation}
    \label{eq:ab2}
    S_k(B):\;
    \begin{tikzcd}
      R_k \arrow[r, rightharpoonup, shift left = 0.3ex, "d_k'"] & \arrow[l, rightharpoonup, shift left = 0.3ex, "c_k'"] V_k[-\lambda_k] \arrow[r, rightharpoonup, shift left = 0.3ex, "\alpha'"] & \arrow[l, rightharpoonup, shift left = 0.3ex, "\beta'"] V_l[\lambda_l+ \lambda_k] \arrow[r, rightharpoonup, shift left = 0.3ex, "c_l"] & \arrow[l, rightharpoonup, shift left = 0.3ex, "d_l"] R_l,
    \end{tikzcd}
  \end{equation}
Note that to determine the maps between $V_l$ and $R_l$ in the above diagram, we have used the condition~\ref{con:c1} in Definition~\ref{def:zk}.

  By the condition~\ref{con:c2} in Definition~\ref{def:zk}, maps in (\ref{eq:ab1}) and (\ref{eq:ab2}) satisfy the short exact sequence,
  \begin{equation}
    \label{eq:sk}
    \begin{tikzcd}[ampersand replacement=\&]
      0 \arrow[r] \& V_k \arrow[r, "{\renewcommand\arraystretch{0.7}\begin{pmatrix}
          c_k' \\ \alpha'
        \end{pmatrix}}"] \& R_k \oplus V_l \arrow[r, "{\renewcommand\arraystretch{0.7}\setlength\arraycolsep{2pt}\begin{pmatrix}
          d_k & \varepsilon\cdot\beta
        \end{pmatrix}}"] \&[1em] V_k \arrow[r] \& 0.
    \end{tikzcd}
  \end{equation}
Moreover, due to our choice of the volume form on $R_k$, the torsion of (\ref{eq:sk}) is $1$.

  \subparagraph{The diagram for $S_lS_k(B)$.}
  Since $(S_k(B),S_lS_k(B))\in Z_l$, by using Remark \ref{rk:bb}~\ref{rk:bb1} again, we have
  \begin{equation}
    \label{eq:ab3}
    S_lS_k(B):\;
    \begin{tikzcd}
      R_k \arrow[r, rightharpoonup, shift left = 0.3ex, "d_k'"] & \arrow[l, rightharpoonup, shift left = 0.3ex, "c_k'"] V_k[\lambda_l] \arrow[r, rightharpoonup, shift left = 0.3ex, "\alpha''"] & \arrow[l, rightharpoonup, shift left = 0.3ex, "\beta''"] V_l[-\lambda_l- \lambda_k] \arrow[r, rightharpoonup, shift left = 0.3ex, "c_l''"] & \arrow[l, rightharpoonup, shift left = 0.3ex, "d_l''"] R_l.
    \end{tikzcd}
  \end{equation}
  By the condition~\ref{con:c2} in Definition~\ref{def:zk}, maps in (\ref{eq:ab2}) and (\ref{eq:ab3}) satisfy the short exact sequence,
\begin{equation}
    \label{eq:skl}
    \begin{tikzcd}[ampersand replacement=\&]
      0 \arrow[r] \& V_l \arrow[r, "{\renewcommand\arraystretch{0.7}\begin{pmatrix}
          c_l'' \\ \beta''
        \end{pmatrix}}"] \& R_l \oplus V_k \arrow[r, "{\renewcommand\arraystretch{0.7}\setlength\arraycolsep{2pt}\begin{pmatrix}
          d_l & -\varepsilon\cdot\alpha'
        \end{pmatrix}}"] \&[2em] V_l \arrow[r] \& 0.
    \end{tikzcd}
  \end{equation}
whose torsion is also $1$.
  Then, by Corollary~\ref{cor:la}, (\ref{eq:sk}) and (\ref{eq:skl}) imply (using $V_k$ to replace $U$ in Corollary~\ref{cor:la}) the short exact sequence
  \begin{equation}
    \label{eq:rkl}
    \begin{tikzcd}[ampersand replacement=\&]
      0 \arrow[r] \& V_l \arrow[r, "{\renewcommand\arraystretch{0.7}\begin{pmatrix}
          c_l'' \\ c'_k\beta''
        \end{pmatrix}}"] \& R_l \oplus R_k \arrow[r, "{\renewcommand\arraystretch{0.7}\setlength\arraycolsep{2pt}\begin{pmatrix}
          \varepsilon\cdot\beta d_l & d_k
        \end{pmatrix}}"] \&[2em] V_k \arrow[r] \& 0,
    \end{tikzcd}
  \end{equation}
  the torsion of which is $(-1)^{\dim{V_k}\dim{V_l} + \dim{V_l}}$.

  \subparagraph{The diagram for $S_kS_lS_k(B)$.}
  Since $(S_kS_lS_k(B), S_lS_k(B))\in Z_k$, by using Remark \ref{rk:bb}~\ref{rk:bb1} again, we have
  \begin{equation}
    \label{eq:ab4}
    S_kS_lS_k(B):\;
    \begin{tikzcd}
      R_k \arrow[r, rightharpoonup, shift left = 0.3ex, "d_k^{\pt}"] & \arrow[l, rightharpoonup, shift left = 0.3ex, "c_k^{\pt}"] V_k[-\lambda_l] \arrow[r, rightharpoonup, shift left = 0.3ex, "\alpha^{\pt}"] & \arrow[l, rightharpoonup, shift left = 0.3ex, "\beta^{\pt}"] V_l[-\lambda_k] \arrow[r, rightharpoonup, shift left = 0.3ex, "c_l''"] & \arrow[l, rightharpoonup, shift left = 0.3ex, "d_l''"] R_l.
    \end{tikzcd}
  \end{equation}
  By the condition~\ref{con:c2} in Definition~\ref{def:zk}, maps in (\ref{eq:ab3}) and (\ref{eq:ab4}) satisfy the short exact sequence,
  \begin{equation}
    \label{eq:sklk}
    \begin{tikzcd}[ampersand replacement=\&]
      0 \arrow[r] \& V_k \arrow[r, "{\renewcommand\arraystretch{0.7}\begin{pmatrix}
          c_k^{\pt} \\ \alpha^{\pt}
        \end{pmatrix}}"] \& R_k \oplus V_l \arrow[r, "{\renewcommand\arraystretch{0.7}\setlength\arraycolsep{2pt}\begin{pmatrix}
          d_k' & \varepsilon\cdot\beta''
        \end{pmatrix}}"] \&[2em] V_k \arrow[r] \& 0,
    \end{tikzcd}
  \end{equation}
  whose torsion is equal to $1$.
By Corollary~\ref{cor:la} again, (\ref{eq:skl}) and (\ref{eq:sklk}) imply (using $V_l$ to replace $U$ in Corollary~\ref{cor:la}) another short exact sequence
  \begin{equation}
    \label{eq:rkl2}
    \begin{tikzcd}[ampersand replacement=\&]
      0 \arrow[r] \& V_k \arrow[r, "{\renewcommand\arraystretch{0.7}\begin{pmatrix}
          c_k^{\pt} \\ c_l''\alpha^{\pt}
        \end{pmatrix}}"] \&[2em] R_k \oplus R_l \arrow[r, "{\renewcommand\arraystretch{0.7}\setlength\arraycolsep{2pt}\begin{pmatrix}
          -\varepsilon\cdot\alpha'd_k' & d_l
        \end{pmatrix}}"] \&[2em] V_l \arrow[r] \& 0,
    \end{tikzcd}
  \end{equation}
  whose torsion is $(-1)^{\dim{V_l}\dim{V_k} + \dim{V_k}}$.

  \subparagraph{The diagram for $S_l(B)$.}
  Since $(B,S_l(B))\in Z_l$, by Remark~\ref{rk:bb}~\ref{rk:bb1}, we have
\begin{equation}
    \label{eq:ab2p}
    S_l(B):\;
    \begin{tikzcd}
      R_k \arrow[r, rightharpoonup, shift left = 0.3ex, "d_k"] & \arrow[l, rightharpoonup, shift left = 0.3ex, "c_k"] V_k[\lambda_k+ \lambda_l] \arrow[r, rightharpoonup, shift left = 0.3ex, "\dot{\alpha}"] & \arrow[l, rightharpoonup, shift left = 0.3ex, "\dot{\beta}"] V_l[-\lambda_l] \arrow[r, rightharpoonup, shift left = 0.3ex, "\dot{c}_l"] & \arrow[l, rightharpoonup, shift left = 0.3ex, "\dot{d}_l"] R_l.
    \end{tikzcd}
  \end{equation}
  By the condition~\ref{con:c2} in Definition~\ref{def:zk}, maps in (\ref{eq:ab1}) and (\ref{eq:ab2p}) satisfy the short exact sequence,
  \begin{equation}
    \label{eq:sl}
    \begin{tikzcd}[ampersand replacement=\&]
      0 \arrow[r] \& V_l \arrow[r, "{\renewcommand\arraystretch{0.7}\begin{pmatrix}
          \dot{c_l} \\ \dot{\beta}
        \end{pmatrix}}"] \& R_l \oplus V_k \arrow[r, "{\renewcommand\arraystretch{0.7}\setlength\arraycolsep{2pt}\begin{pmatrix}
          d_l & -\varepsilon\cdot\alpha
        \end{pmatrix}}"] \&[2em] V_l \arrow[r] \& 0,
    \end{tikzcd}
  \end{equation}
  whose torsion is also $1$.

\paragraph{Step 2.} So far, we have introduced quite a few maps.
  In this step, based on these maps, we will construct a nice candidate for ``$S_kS_l(B)$'', that is, we will find $A\in \mm^{-1}(0)$ such that $(S_l(B),A)\in Z_k$.

  To find such a point $A$, we use the following result due to Maffei.
  \begin{lemma}[{\cite[p.680]{Maffei:2002aa}}]
    \label{lm:maf}
    There exist unique linear maps $x:V_k\rightarrow V_l$ and $y:V_l\rightarrow V_k$ such that the following equations hold,
    \begin{equation*}
      \left\{
        \begin{aligned}
          \dot{c}_lx &= c_l''\alpha^{\pt} \\
          \dot{\beta} x &= -\varepsilon\cdot d_kc_k^{\pt}
        \end{aligned}\right.\quad\text{and}\quad
      \left\{
        \begin{aligned}
          xy & = \dot{\alpha}\dot{\beta} - \varepsilon\cdot(\lambda_k + \lambda_l)\id_{V_l} \\
          c_k^{\pt}y & = c_k\dot{\beta}
        \end{aligned}\right. .
    \end{equation*}
  \end{lemma}

  By using the maps $x$ and $y$ in Lemma~\ref{lm:maf}, we choose $A\in \bfM$, which corresponds to the following diagram,
  \begin{equation}
    \label{eq:ab3p}
    A:\;
    \begin{tikzcd}
      R_k \arrow[r, rightharpoonup, shift left = 0.3ex, "d_k^{\pt}"] & \arrow[l, rightharpoonup, shift left = 0.3ex, "c_k^{\pt}"] V_k[-\lambda_k- \lambda_l] \arrow[r, rightharpoonup, shift left = 0.3ex, "x"] & \arrow[l, rightharpoonup, shift left = 0.3ex, "y"] V_l[\lambda_k] \arrow[r, rightharpoonup, shift left = 0.3ex, "\dot{c}_l"] & \arrow[l, rightharpoonup, shift left = 0.3ex, "\dot{d}_l"] R_l.
    \end{tikzcd}
  \end{equation}
  More precisely, for $h\in \DQw$ such that $\hd{h} \neq k,l$ and $\tl{h} \ne k,l$, we define $A_h = B_h$.
  For the remaining components of $A$, they are defined according to diagram (\ref{eq:ab3p}).
  We are going to show that $A\in \mm^{-1}(0)$ and $(S_l(B), A)\in Z_k$.

  In fact, since $S_l(B) \in \bfN_{kl}$, if we can show that the pair $(S_l(B), A)$ satisfies all three conditions in Definition~\ref{def:zk}, Remark~\ref{rk:bb} \ref{rk:bb2} implies $A\in \bfN_{kl}\subseteq \mm^{-1}(0)$ automatically.
  Hence, in the following, we only need to prove that $(S_l(B), A)$ satisfies all three conditions in Definition~\ref{def:zk}.

  Firstly, by our definition of $A$, $(S_l(B),A)$ satisfies the condition~\ref{con:c1} trivially.
  Next, we check that this pair also satisfies the condition~\ref{con:c2}, that is, we need to show that the following chain complex is a torsion $1$ short exact sequence,
  \begin{equation}
    \label{eq:slk}
    \begin{tikzcd}[ampersand replacement=\&]
      0 \arrow[r] \& V_k \arrow[r, "{\renewcommand\arraystretch{0.7}\begin{pmatrix}
          c_k^{\pt} \\ x
        \end{pmatrix}}"] \& R_k \oplus V_l \arrow[r, "{\renewcommand\arraystretch{0.7}\setlength\arraycolsep{2pt}\begin{pmatrix}
          d_k & \varepsilon\cdot \dot{\beta}
        \end{pmatrix}}"] \&[1em] V_k \arrow[r] \& 0.
    \end{tikzcd}
  \end{equation}
The exactness of (\ref{eq:slk}) follows from Lemma~\ref{lm:maf} immediately.
  We only need to calculate its torsion.

  Recall that $(B,S_k(B))\in Z_k$.
  Due to the condition~\ref{con:c3} in Definition~\ref{def:zk}, the maps in diagrams (\ref{eq:ab1}) and (\ref{eq:ab2}) satisfy
  \begin{equation}
    \label{eq:c3-ab12}
    \alpha'd'_k = \alpha d_k.
  \end{equation}
  Now by applying Corollary~\ref{cor:la} to (\ref{eq:sl}) and (\ref{eq:slk}) (using $V_l$ to replace $U$ in Corollary~\ref{cor:la}), due to (\ref{eq:c3-ab12}), the output short exact sequence is (\ref{eq:rkl2}) exactly.
  Therefore, by using the torsion of (\ref{eq:rkl2}) and (\ref{eq:sl}), we conclude that the torsion of (\ref{eq:slk}) is $1$ as expected.

  What is left is to check $(S_l(B),A)$ satisfying the condition~\ref{con:c3}.
  In view of Lemma~\ref{lm:maf}, we only need to show the following two equalities,
  \begin{subequations}
    \begin{align}
      c^{\pt}_kd^{\pt}_k &= c_kd_k - (\lambda_k + \lambda_l)\id_{R_k}, \label{eq:c31}\\
      xd^{\pt}_k &= \dot{\alpha}d_k.\label{eq:c32}
    \end{align}
  \end{subequations}

  We first show (\ref{eq:c31}).
  Since $(B,S_k(B))\in Z_k$, the condition~\ref{con:c3} implies the following equality between maps in (\ref{eq:ab1}) and (\ref{eq:ab2}),
  \begin{equation}
    \label{eq:c31-1}
    \begin{pmatrix}
      c'_k \\ \alpha'
    \end{pmatrix}
    \begin{pmatrix}
      d'_k & \beta'
    \end{pmatrix}
    =
    \begin{pmatrix}
      c_k \\ \alpha
    \end{pmatrix}
    \begin{pmatrix}
      d_k & \beta
    \end{pmatrix}
    - \lambda_k\id_{R_k \oplus V_k}.
  \end{equation}
  For the same reason, using $(S_lS_k(B),S_kS_lS_k(B))\in Z_k$, we have
  \begin{equation}
    \label{eq:c31-2}
    \begin{pmatrix}
      c^{\pt}_k \\ \alpha^{\pt}
    \end{pmatrix}
    \begin{pmatrix}
      d^{\pt}_k & \beta^{\pt}
    \end{pmatrix}
    =
    \begin{pmatrix}
      c'_k \\ \alpha''
    \end{pmatrix}
    \begin{pmatrix}
      d'_k & \beta''
    \end{pmatrix}
    - \lambda_l\id_{R_k \oplus V_k}.
  \end{equation}
  (\ref{eq:c31}) follows from (\ref{eq:c31-1}) and (\ref{eq:c31-2}) immediately.

To show (\ref{eq:c32}), we first note that (\ref{eq:c31-2}) implies the following analog of (\ref{eq:c3-ab12}),
  \begin{equation}
    \label{eq:c32-1}
    \alpha^{\pt}d_k^{\pt} =  \alpha''d_k'.
  \end{equation}
  Similarly, since $(S_k(B), S_lS_k(B))\in Z_l$, $(B,S_l(B))\in Z_l$, by using the condition~\ref{con:c3} in Definition~\ref{def:zk}, we have
  \begin{align}
    c_l \alpha' & = c_l'' \alpha'',\label{eq:c32-2}\\
    c_l \alpha  & = \dot{c}_l\dot{\alpha}.\label{eq:c32-3}
  \end{align}

  Meanwhile, from the diagram (\ref{eq:ab2p}), we know that $\mu_k(S_l(B)) = (\lambda_k + \lambda_l)\id_{V_k}$.
  By the definition of $\mu$, this equality can be written as
  \begin{equation}
    \label{eq:c32-4}
    (\lambda_k + \lambda_l)\id_{V_k} = \mu_k(S_l(B)) = d_kc_k + \varepsilon\cdot\dot{\beta}\dot{\alpha}.
  \end{equation}

  Since we have shown (\ref{eq:c31}), combining with Lemma~\ref{lm:maf} and (\ref{eq:c3-ab12}), (\ref{eq:c32-1}-\ref{eq:c32-4}), we calculate as follows,
  \begin{equation*}
    \dot{c}_lxd^{\pt}_k \overset{\text{lm.~\ref{lm:maf}}}{\longeq} c_l''\alpha^{\pt}d^{\pt}_k
    \overset{\text{(\ref{eq:c32-1})}}{\longeq} c''_l\alpha'' d'_k
    \overset{\text{(\ref{eq:c32-2})}}{\longeq} c_l\alpha' d'_k
    \overset{\text{(\ref{eq:c3-ab12})}}{\longeq} c_l\alpha d_k \overset{\text{(\ref{eq:c32-3})}}{\longeq} \dot{c}_l\dot{\alpha}d_k,
  \end{equation*}
  and
  \begin{equation*}
\dot{\beta}xd^{\pt}_k \overset{\text{lm.~\ref{lm:maf}}}{\longeq} -\varepsilon\cdot d_kc^{\pt}_kd^{\pt}_k \overset{(\ref{eq:c31})}{\longeq} -\varepsilon\cdot d_k(c_kd_k - (\lambda_k + \lambda_l)\id_{R_k})  \overset{(\ref{eq:c32-4})}{\longeq} \dot{\beta}\dot{\alpha}d_k.
  \end{equation*}
The above two equalities imply the following equality,
  \begin{equation}
    \label{eq:c32-5}
    \begin{pmatrix}
      \dot{c}_l \\ \dot{\beta}
    \end{pmatrix}xd^{\pt}_k =
    \begin{pmatrix}
      \dot{c}_l \\ \dot{\beta}
    \end{pmatrix}\dot{\alpha}d_k.
  \end{equation}
  Since (\ref{eq:sl}) is exact, $\tsp{(\dot{c}_l,\dot{\beta})}$ must be injective.
  Hence, (\ref{eq:c32}) follows from (\ref{eq:c32-5}).

  \paragraph{Step 3.}
  The last step for the proof of $Z_{klk}\subseteq Z_{lkl}$ is to check that $(A,S_kS_lS_k(B)) \in Z_l$.
  However, this can be done by using essentially the same argument that we have used to prove $(S_l(B),A)\in Z_k$.

  As before, the upshot is to show the following chain complex is exact and has torsion $1$ property,
  \begin{equation}
    \label{eq:slkl}
    \begin{tikzcd}[ampersand replacement=\&]
      0 \arrow[r] \& V_l \arrow[r, "{\renewcommand\arraystretch{0.7}\begin{pmatrix}
          c_l'' \\ \beta^{\pt}
        \end{pmatrix}}"] \& R_l \oplus V_k \arrow[r, "{\renewcommand\arraystretch{0.7}\setlength\arraycolsep{2pt}\begin{pmatrix}
          \dot{d}_l & -\varepsilon\cdot x
        \end{pmatrix}}"] \&[2em] V_l \arrow[r] \& 0.
    \end{tikzcd}
  \end{equation}
  Recall that, in the proof of $(S_l(B),A)\in Z_k$, we have shown that (\ref{eq:slk}) has these two properties.
  To show that (\ref{eq:slkl}) also has these two properties, we can use a similar method.
  The only difference is that in current situation, (\ref{eq:rkl}) plays the role that is played by (\ref{eq:rkl2}) before.
  The details are left to readers.

  The proof of Proposition~\ref{prop:klk} is finished.
\end{proof}

\vspace{3ex}
Now we have verified $\aset{S_k}_{k\in I}$ satisfying all the Coxeter relations and the proof of Theorem~\ref{thm:main-re} is completed.
\begin{remark}
  If we check the whole proof for the existence of the Weyl group action on $\cM(\bfv,\bfw)$, as readers may have noticed, the dimension vectors condition $\bfw = \bfC \bfv$ is not always relevant.
  It seems possible to relax this condition and prove a result more similar to Theorem~\ref{thm:nweyl}.
\end{remark}

\section{The comparison of two Weyl group actions: a case study}
\label{sec:comparison}

\subsection{}
This section consists of two parts.
In the first part, \S~\ref{sub:goveru} and \S~\ref{sub:weyl-gu}, we review another Weyl group action, i.e.\ the \grk\ action on $\afbc$ (or $\big(\rT^*(G^{\mathrm{ad}}/U)\big)_{\aff}$) mentioned in Introduction.
In the second part, \S~\ref{sub:an-type} and \S~\ref{sub:pf-thm2}, we concentrate on a special case, $G = \SL{n+1}$.
We recall a result in~\cite{Dancer:2013aa} which asserts that for this special $G$, $\afbc$ is also a DKS type variety.
Then for this special example of DKS type varieties, we compare the Weyl group action constructed in the previous section with the \grk\ action.
More specifically, we prove Theorem~\ref{thm:main2}.

\subsection{The geometry of $G/U$.}
\label{sub:goveru}
We review some geometric properties of the basic affine variety to be used later.
Especially, we recall a canonical construction of $G/U$ and $\rT^*(G/U)$ after Ginzburg and Kazhdan.
To avoid too many details, we only formulate the construction assuming $G$ to be a connected semi-simple group with a trivial center, i.e.\ $G$ is of adjoint type. 
But we would like to point out that, by proceeding as in~\cite[\S~5.4]{Ginzburg:2018la}, it is also possible to give a similar construction for the ``simply-connected'' case.

Let $\kg = \mathrm{Lie}(G)$ and $\cB$ be the flag variety of Borel subgroups of $B\subseteq G$, or equivalently, Borel subalgebras of $\kb\subseteq \kg$.
Given a Borel subalgebra $\kb\in \cB$, we always use $B$ to denote the corresponding Borel subgroup.
The tori $B/[B,B]$ for any $B\in \cB$ are canonically isomorphic, which is denoted by $T$ and called the \emph{universal Cartan torus} for $G$.
Let $\kt = \mathrm{Lie}(T)$ and $\Lambda^* \allowbreak= \Hom{T}{\bC^*}\subseteq \kt^*$ be the root lattice of $\kg$.
For any Cartan subgroup $H\subseteq G$,
\begin{equation*}
  W_H \coloneqq \opN_G(H)/H
\end{equation*}
is defined to be the Weyl group associated with $H$.
Recall that there is another Weyl group $\bW$ associated with the Dynkin diagram of $G$, which we call the \emph{universal Weyl group}.
By choosing a Borel subgroup $B \supseteq H$, $W_H$ is canonically isomorphic to $\bW$.
Therefore, there is a canonical $\bW$ action on $T$.

Given a Borel subalgebra $\kb$, let $\ku(\kb) \coloneqq [\kb,\kb]$ be the nilradical of $\kb$.
Then we define
\begin{equation*}
  \ka(\kb) \coloneqq \ku(\kb)/[\ku(\kb),\ku(\kb)].
\end{equation*}
The weights of the $T$ action on $\ka(\kb)$ are called simple roots.
Denote simple roots by $\aset{\lambda_i}_{i\in I}$, where $I$ is the vertex set of the Dynkin diagram.
Define a partial order on $\Lambda^*$ as follows: for $\mu_1,\mu_2\in \Lambda^*$, $\mu_1 \ge \mu_2$ if and only if $\mu_1- \mu_2 = \sum_{i\in I} \bZ_{\ge 0}\cdot \lambda_i$.
For a fixed $\kb\in \cB$, there exists a canonical $\kb$ stable filtration $\kg^{\ge \mu,\kb}$, $\mu \in \Lambda^*$, of $\kg$, such that $\gr^{\mu,\kb} \kg \coloneqq \kg^{\ge \mu,\kb}/ \kg^{> \mu,\kb}$ is a representation of the weight $\mu$ under the natural action of $T = B/[B,B]$.
For example, taking $\mu=0$, one can check that $\kg^{\ge 0,\kb} = \kb$ and $\gr^{0,\kb}\kg = \kt$.
Using this filtration, we define a special subspace of $\kg$,
\begin{equation*}
  \kd(\kb) \coloneqq \sum_{i\in I} \kg^{\ge -\lambda_i,\kb}.
\end{equation*}
In a sense that we will clarify a moment later, $\kd(\kb)/\kb$ can be seen as a dual to $\ka(\kb)$.
Like $\ka(\kb)$, $T$ also acts on $\kd(\kb)/\kb$ naturally.
The unique open dense orbit of the $T$ action on $\ka(\kb)$ (resp.\ $\kd(\kb)/\kb$) is denoted by $\bO(\kb)$ (resp.\ $\bO_{-}(\kb)$).
We will abuse notations a little by writing the preimage of $\bO_-(\kb)$ in $\kd(\kb)$ as $\bO_-(\kb) + \kb$.

Like $\cB$, we can construct two varieties consisting of a family of $\bO(\kb)$ or $\bO_{-}(\kb)$,
\begin{equation*}
  \widetilde{\cB} \coloneqq \aset{(\kb,\sfs) | \kb\in \cB,\, \sfs\in \bO(\kb)},\quad \widetilde{\cB}_- \coloneqq \aset{(\kb,\sfs) | \kb\in \cB,\, \sfs\in \bO_-(\kb)}.
\end{equation*}
$G$ acts on $\cB$, $\widetilde{\cB}$, $\widetilde{\cB}_-$ by conjugation.
And the $T$ action on $\bO(\kb)$ (resp.\ $\bO_-(\kb)$) induces a $T$ action on $\widetilde{\cB}$ (resp.\ $\widetilde{\cB}_-$), which commutes with the $G$ action on it.
Therefore, $\widetilde{\cB}$ (resp.\ $\widetilde{\cB}_-$) is a $T$-principal bundle over $\cB$ and the canonical projection of this bundle is $G$-equivariant.

\begin{example*}
  \label{ex:sl}
  For $G= \PGL{n+1}$, take $\kb\in \cB$ to be the Borel subalgebra consisting of upper triangular matrices.
  The preimages of $\bO(\kb)$ and $\bO_-(\kb)$ in $\ksl(n+1,\bC)$ are matrices of the following forms.
  \begin{equation*}
\text{For }\bO(\kb),\,
\begin{pmatrix}
      0      & e_1    & *      & \cdots       & * \\
      \vdots      & \ddots      &    \ddots  &   \ddots   & \vdots \\
             &         &  \ddots     &   \ddots    & * \\
      \vdots      &        &        &  \ddots      & e_{n} \\
      0      &        & \cdots      &         & 0
    \end{pmatrix},    
    \text{ and for }\bO_-{(\kb)},\,
\begin{pmatrix}
      *      &        & \cdots      &        & * \\
      f_1    &   \ddots    &        &        & \vdots \\
      0      &   \ddots   &   \ddots    &        &    \\
      \vdots     &   \ddots   &   \ddots    & \ddots      &  \vdots\\
      0      & \cdots      & 0      & f_{n}   & *
    \end{pmatrix},
\end{equation*}
  where $e_i$ and $f_i$ are non-zero numbers.
  Taking $t= \diag{z_0,z_1,\cdots,z_{n}}\in \kb$ as an element in $T$, under the adjoint action of $t$ on the above matrices, the entry $e_i$ in the first matrix is scaled by $z_{i}^{-1}z_{i-1}$ and the entry $f_i$ in the second matrix is scaled by $z_{i-1}^{-1}z_{i}$.
\end{example*}

Using the Killing form of $\kg$, we will identify $\kg$ with $\kg^*$.
Since the isotropic subgroup at $(\kb,\sfs)\in \widetilde{\cB}$ is the maximal unipotent subgroup associated with the Lie algebra $\ku(\kb)$, the cotangent bundle of $\widetilde{\cB}$ is\footnote{Here and in the following, we often only define a variety (or a map between two varieties) as a set.
  But all such definitions can be enhanced scheme-theoretically by using suitable base changes.
  For details, one can see~\cite{Ginzburg:2018la}.}
\begin{equation*}
  \rT^*\widetilde{\cB} = \aset{(\kb,\sfs,x)| \kb\in \cB,\,\sfs\in \bO(\kb),\,x\in \ku(\kb)^{\perp}\simeq \kb}.
\end{equation*}
Fixing a choice of base point $(\kb,\sfs)\in \widetilde{\cB}$, thus also fixing a Borel subgroup $B$, the following two $G$-equivariant isomorphisms hold,
\begin{equation}
  \label{eq:2mor}
  \begin{aligned}
    G/U &\simarrow \widetilde{\cB} & G\times_{U} \kb &\simarrow \rT^*\widetilde{\cB} \\
    [g] &\mapsto \big(\Ad{g}{\kb},\Ad{g}{\sfs}\big), & [g,x] &\mapsto \big(\Ad{g}{\kb},\Ad{g}{\sfs},\Ad{g}{x}\big),
  \end{aligned}
\end{equation}
where $U= [B,B]$ and $\Ad{g}{\sfs}\in \bO(\Ad{g}{\kb})$.

Using the transitive $G$ action on $G/U$, $\rT^*(G/U)$ is isomorphic to $G\times_U (\kg/\ku(\kb))^*$.
And by using the Killing form, we have $G\times_U (\kg/\ku(\kb))^* \simeq G\times_{U} \kb$.
Therefore, from now on, we will identify $\rT^*(G/U)$ with $G\times_{U} \kb$.
By using the second isomorphism in (\ref{eq:2mor}), we have an isomorphism
\begin{equation*}
  \Theta_{\sfs}: \rT^*(G/U)\, (\simeq G\times_{U} \kb) \simarrow \rT^*\widetilde{\cB}.
\end{equation*}
But we would like to remind readers that the isomorphism $\Theta_{\sfs}$ obtained in this way \uline{does depend on} the choice of $\sfs$.
Considering this fact, it is reasonable to view $\widetilde{\cB}$ and $\rT^*\widetilde{\cB}$ as canonical models of $G/U$ and $\rT^*(G/U)$ respectively.
\begin{remark}
  \label{rk:sc}
  If $G$ is a simply connected group, one can see that the two morphisms in (\ref{eq:2mor}), as well as $\Theta_{\sfs}$, are finite covering maps, whose fibers are homeomorphic to the center of $G$.
  By using a variation of $\widetilde{\cB}$, denoted by $\widetilde{\cB}^{\mathrm{sc}}$, it is possible to construct isomorphisms like (\ref{eq:2mor}) for a simply connected group.
  Roughly speaking, one needs to choose the fundamental weights $\aset{\varpi_i}$ of $T$ and the corresponding representations $\aset{V_{\varpi_i}}$.
  Then $\widetilde{\cB}^{\mathrm{sc}}$ is defined in the same way as $\widetilde{\cB}$ except that $\sfs$ is replaced by $(\rv_1,\cdots,\rv_n)$, where $\rv_i$ is a highest weight vector of the representation $V_{\varpi_i}$.
  Readers can find more details in~\cite[\S~5.4]{Ginzburg:2018la}.
  We only would like to point out that the isomorphism between $G/U$ and $\widetilde{\cB}^{\mathrm{sc}}$ is actually the same as the isomorphism constructed in~\cite[\S~6]{Guillemin:2002sx}.
\end{remark}

Recall that the Grothendieck-Springer resolution $\widetilde{\kg}\coloneqq\aset{(\kb,x)\in \cB\times \kg| x\in\kb }$ makes the following diagram commute,
\begin{equation*}
  \begin{tikzcd}
    & \arrow[ld, "x\mapsfrom (\kb{,}x) :\pi"'] \widetilde{\kg} \arrow[rd, "\nu:(\kb{,}x)\mapsto x\bmod{\ku(\kb)}"] & \\
    \kg \arrow[rd]  & & \arrow[ld] \kt \\
    &  \kg \dsl G &
  \end{tikzcd},
\end{equation*}
where the GIT quotient $\kg\dsl G$ is isomorphic to $\kt/\bW$ by the Chevalley isomorphism.
The definition of $\rT^*\widetilde{\cB}$ gives a natural projection \begin{alignat*}{2}
  \opp_{\widetilde{\kg}}:& &\rT^*\widetilde{\cB} &\rightarrow \widetilde{\kg}\\
  & & (\kb,\sfs,x) &\mapsto (\kb, x).
\end{alignat*}
Since $G \times T$ acts on $\widetilde{\cB}$, there exists a Hamiltonian $G\times T$ action on $\rT^*\widetilde{\cB}$ with the moment map $\mu_G\times \mu_T$.
Using $\opp_{\widetilde{\kg}}$ and the Killing form, this moment map has the following expression: \footnote{Compared to~\cite{Ginzburg:2018la}, the moment maps here have an extra minus sign because~\cite{Ginzburg:2018la} and this paper use different sign conventions for moment maps.}
\begin{align*}
  \mu_G &= -\pi \circ \opp_{\widetilde{\kg}}: \rT^*\widetilde{\cB} \rightarrow \kg,\\
  \mu_T &= -\nu \circ \opp_{\widetilde{\kg}}:\rT^*\widetilde{\cB} \rightarrow \kt.
\end{align*}

\subsection{The Weyl group action on $(\rT^*\widetilde{\cB})_{\aff}$.}
\label{sub:weyl-gu}

To recall Ginzburg and Kazhdan's construction of the Weyl group action, i.e.\ the \grk\ action, on $(\rT^*\widetilde{\cB})_{\aff}$, we begin with the definition of a special variety $\cZ_{\ras}$.
As in the previous subsection, we also assume that $G$ is of adjoint type.

Let $\kg_{\ras}\subseteq \kg$ be subset of the regular and semi-simple elements.
Then $\rT^*\widetilde{\cB}_{\ras} \coloneqq \mu_G^{-1}(\kg_{\ras})$ (resp.\ $\widetilde{\kg}_{\ras} \coloneqq \pi^{-1}(\kg_{\ras})$) is a Zariski open and dense subvariety of $\rT^*\widetilde{\cB}$ (resp.\ $\widetilde{\kg}$).
Recall that for a pair of Borel subalgebras $\kb,\bkb\in \cB$, they are called \emph{in opposite position} if $\kb \cap \bkb$ is a Cartan subalgebra of $\kg$.
After~\cite[p.10,~Definition]{Ginzburg:2018la}, we call the following incidence variety the \emph{Miura variety} for regular and semi-simple elements.
\begin{equation*}
  \cZ_{\ras} \coloneqq \aset{(\bkb,\kb,x) \in \cB \times \cB \times \kg_{\ras} | x\in \kb \cap [\bO_-(\bkb) + \bkb]}.
\end{equation*}
\begin{remark*}
  By~\cite[Proposition~4.1.1]{Ginzburg:2018la}, if a point $(\bkb,\kb,x)$ lies in $\cZ_{\ras}$, the pair $(\bkb,\kb)$ must be in opposite position.
  Moreover, in~\cite{Ginzburg:2018la}, the Miura variety is defined for regular elements, not necessarily semi-simple.
  But for our applications, the above smaller variety is sufficient.
\end{remark*}

For a pair of Borel subalgebras $(\kb,\bkb)$ in opposite position, we have a triangular decomposition $\kg = \ku(\kb) \oplus (\kb \cap \bkb) \oplus \ku(\bkb)$ and a diagram
\begin{equation*}
  \kt = \bkb/\ku(\bkb) \twoheadleftarrow  \bkb \hookleftarrow \bkb \cap \kb \hookrightarrow \kb \twoheadrightarrow \kb/\ku(\kb) = \kt.
\end{equation*}
Using the above diagram, there is an isomorphism between the leftmost and the rightmost copy of $\kt$: $t\mapsto w_0(t)$, where $w_0$ is the longest element in the Weyl group.
Similarly, the following two maps are also isomorphisms,
\begin{equation*}
  \begin{tikzcd}
    \kd(\bkb)/\bkb & \arrow[l, "\kappa_-"', "\simeq"] \kd(\bkb) \cap \ku(\kb) \arrow[r, "\kappa_+", "\simeq"'] & \ka(\kb).
  \end{tikzcd}
\end{equation*}
Compositing the above two isomorphisms, we have a new isomorphism $\kappa_{\bkb,\kb}$,
\begin{equation*}
  \kappa_{\bkb,\kb} \coloneqq \kappa_+\kappa_-^{-1}:\kd(\bkb)/\bkb \simarrow \ka(\kb).\end{equation*}
$\kappa_{\bkb,\kb}$ intertwines the $T$ action in the following way,
\begin{equation*}
  {\kappa_{\bkb,\kb}(t\mathsf{f}) = w_0^{-1}(t)\kappa_{\bkb,\kb}(\mathsf{f}),\quad \text{if } t\in T,\mathsf{f}\in \kd(\bkb)/\bkb}.
\end{equation*}
Hence, $\kappa_{\bkb,\kb}$ induces an isomorphism from $\bO_-{(\bkb)}$ to $\bO(\kb)$, which is denoted by the same symbol $\kappa_{\bkb,\kb}$.

For any $(\bkb,\kb,x)\in \cZ_{\ras}$, the isotropic subgroup at $x$, $H_x$, is a Cartan torus of $G$, which is contained in $B$.
Therefore, $H_x$ is isomorphic to $B/[B,B]$ and isomorphic to $T$ consequently.
As a result, $T$ can act on $\cZ_{\ras}$ in the following way.
For any $t\in T$, let $h \in H_x$ be the element corresponding to $t$.
\begin{align*}
  T \times \cZ_{\ras} & \rightarrow \cZ_{\ras} \\
  t\times (\bkb,\kb,x) & \mapsto t\cdot(\bkb,\kb,x) \coloneqq (\Ad{h}{\bkb},\kb,x).
\end{align*}

Using $\kappa_{\bkb,\kb}$, we can identify $\cZ_{\ras}$ with $\rT^*\widetilde{\cB}_{\ras}$ as follows.
\begin{proposition}[{\cite[Corollary~5.2.4]{Ginzburg:2018la}}]
  \label{prop:gk}
  The following map is a $T$-equivariant isomorphism,
  \begin{alignat*}{2}
    \kappa:& & \cZ_{\ras} &\rightarrow \rT^*\widetilde{\cB}_{\ras},\\
    & & (\bkb,\kb,x) &\mapsto (\kb, \kappa_{\bkb,\kb}(x \bmod \bkb), x).
  \end{alignat*}
\end{proposition}
By the definition of $\kappa$, $\opp_{\widetilde{\kg}}\circ \kappa$ coincides with the natural projection $(\bkb,\kb,x)\mapsto (\kb,x)$.
Without any ambiguities, we also use the symbol $\opp_{\widetilde{\kg}}$ for this composition map.

By Proposition~\ref{prop:gk}, to define a Weyl group action on $\rT^*\widetilde{\cB}_{\ras}$, we only need to construct a Weyl group action on $\cZ_{\ras}$.
As before, for any $(\bkb,\kb,x) \in \cZ_{\ras}$, let $H_x\subseteq B$ be the isotropic subgroup at $x$.
For any $n^{w} \in \opN_G(H_x)$ corresponding to an element $w\in \bW \simeq W_{H_x}$ via the pair $(H_x,B)$, let
\begin{equation*}
  \kb^{w} \coloneqq \Ad{n^{w}}{\kb}.
\end{equation*}
Then as in~\cite[\S~4.2]{Ginzburg:2018la}, we define the image of $(\bkb,\kb, x)$ under the action of $w$ to be
\begin{equation}
  \label{eq:defw1}
  w (\bkb, \kb, x) \coloneqq (\bkb,\kb^{w^{-1}},x).\footnote{The Weyl group action given here is a little different from~\cite{Ginzburg:2018la} because we want to ensure the Weyl group action to be a left action.}
\end{equation}
By definition, the $G$ action and the $\bW$ action on $\cZ_{\ras}$ commute.
Moreover, for any $t\in T$, we have $w \big(t\cdot(\bkb, \kb, x)\big) = w(t)\cdot \big(w (\bkb, \kb, x)\big)$. In other words, there is a $\bW \ltimes T$ action on $\cZ_{\ras}$.

Recall that the Weyl group action on $\widetilde{\kg}_{\ras}$ is defined in a similar way, that is, $(\kb,x) \mapsto (\kb^{w^{-1}},x)$ for any $(\kb,x) \in \widetilde{\kg}_{\ras}$.
Therefore, the natural projection $\opp_{\widetilde{\kg}}$ is equivariant with respect to the Weyl group actions on $\cZ_{\ras}$ and $\widetilde{\kg}_{\ras}$.
Especially, this implies the $G$ moment map $\mu_G: \rT^*\widetilde{\cB}_{\ras} \rightarrow \kg_{\ras}$ is $\bW$ invariant (the $\bW$ action on $\kg$ is trivial).
In the same way, we can see that the $T$ moment map $\mu_T: \rT^*\widetilde{\cB}_{\ras} \rightarrow \kt_{\rr}$ is $\bW$-equivariant.

A more or less annoying fact is that the Weyl group action on $\rT^*\widetilde{\cB}_{\ras}$ can not be extended to a group action on $\rT^*\widetilde{\cB}$.
Nevertheless, by~\cite[Theorem~5.2.7]{Ginzburg:2018la}, it is possible to extend the above Weyl group action uniquely to a group action on a larger variety $(\rT^*\widetilde{\cB})_{\aff} \coloneqq \Spec{\bC[\rT^*\widetilde{\cB}]}$, of which $\rT^*\widetilde{\cB}_{\ras}$ is still an open dense subvariety.
In~\cite{Ginzburg:2018la}, Ginzburg and Kazhdan refer the Weyl group action on $(\rT^*\widetilde{\cB})_{\aff}$ obtained in this way as the \emph{quasi-classical Gelfand-Graev action}.
But, as we have done in Introduction, we will simply call this action the \grk\ action for short.

As Remark~\ref{rk:sk}, the construction in this subsection can also be extended to the simply-connected case, c.f.~\cite[\S~5.4]{Ginzburg:2018la}.

\subsection{}
\label{sub:an-type}
In the rest of this paper, we take $G = \SL{n+1}$, $G^{\mathrm{ad}}= \PGL{n+1}$ and $\kb$ to be the Borel subalgebra of upper triangular matrices.
Therefore, $U$ is the group of invertible upper triangular matrices whose diagonal entries are equal to $1$.

In this subsection, we recall an isomorphism $\Xi$ from a DKS type variety of $\rA_{n}$ type to $\afbc$ following~\cite[\S~7]{Dancer:2013aa}.
Then, via the isomorphism $\Xi$, we will calculate the Weyl group action appeared in Theorem~\ref{thm:main}, or Theorem~\ref{thm:main-re}, explicitly.

\subsubsection{The isomorphism $\Xi$.}
For an $\rA_n$ quiver $Q$, let $\bfv = \tsp{(n,n-1,\cdots,1)}$, $\bfw = \tsp{(n+1,0,\cdots,0)}$ be dimension vectors.\footnote{We choose to write the components of $\bfv$ in a descending order, which is in conformity with the convention used in~\cite{Nakajima:1994aa} but different from the convention in~\cite{Dancer:2013aa}.}
One can check that $\bfw = \bfC\bfv$.
We collect all the data in the following diagram for $\DQ$, Figure~\ref{fig:g-over-u}, in which we choose the orientation of arrows associated with $\beta_k$, $1\le k \le n-1$, to be positive, i.e.\ $\beta_k\in Q$.
Moreover, we assume that the vector space $\bC^l$ associated with the vertex $l$, $1\le l \le n+1$, has the standard volume form.
\begin{figure}[htbp]
  \centering
  \begin{equation*}
    \begin{tikzcd}
      {\scriptstyle \bC^n} &[-0.8em] {\scriptstyle \bC^{n-1}} &[-0.6em] &[-0.3em] {\scriptstyle \bC^{2}} &[-0.4em] {\scriptstyle \bC^1} \\[-2em]
      {\bullet} \arrow[r, rightharpoonup, dashed, shift left = 0.3ex, "\beta_1"] \arrow[d, rightharpoonup, shift left = 0.3ex, "\alpha_0"] & \arrow[l, rightharpoonup, shift left = 0.3ex, "\alpha_1"] \bullet \arrow[r, rightharpoonup, dashed, shift left = 0.3ex] & \arrow[l, rightharpoonup, shift left = 0.3ex, ] \cdots \arrow[r, rightharpoonup, dashed, shift left = 0.3ex] & \arrow[l, rightharpoonup, shift left = 0.3ex] \bullet \arrow[r, rightharpoonup, dashed, shift left = 0.3ex, "\beta_{n-1}"] & \arrow[l, rightharpoonup, shift left = 0.3ex, "\alpha_{n-1}"] \bullet  \\[0.5em]
      \circ \arrow[u, rightharpoonup, shift left = 0.3ex, "\beta_0"]  \arrow[r, phantom, "{\scriptstyle \bC^{n+1}}" near start] & \phantom{a} & & &
    \end{tikzcd}.
  \end{equation*}
  \caption{A quiver associated with $\rT^*(G/U)$.}
  \label{fig:g-over-u}
\end{figure}

With the above settings, the locus $\mm^{-1}(0) \subseteq \bfM$ satisfies the following equations,
\begin{equation}
  \label{eq:slmm}
  - \alpha_k\beta_k + \beta_{k-1}\alpha_{k-1} = \tau_kI_{n-k+1},\quad\text{for }1\le k \le n,
\end{equation}
where $\alpha_n = \beta_n = 0$, $\tau_k\in \bC$ and $I_{n-k+1}$ is an identity matrix of size $n -k +1$.
Consider an open and dense subvariety of $\mm^{-1}(0)$,
\begin{multline*}
  \mm^{-1}_{\surj}(0) \coloneqq \{(\alpha_0,\beta_0,\cdots,\alpha_{n-1},\beta_{n-1}) \in \mm^{-1}(0)|\beta_k\text{ is surjective} \\
  \text{for any }0 \le k \le n-1 \}.
\end{multline*}

Now we choose some matrices of special forms, $\beta^0_k = (\zm{n-k}{1}\, |\, I_{n-k})$, $0\le k \le n-1$.
For any point $(\alpha_0,\beta_0,\cdots,\alpha_{n-1},\beta_{n-1})\in \mm^{-1}_{\surj}(0)$, we say a family of invertible matrices
\begin{equation*}
  g_k\in \SL{n+1-k},\; 0 \le k \le n,
\end{equation*}
is \emph{fitting} for $(\alpha_0,\beta_0,\cdots,\alpha_{n-1},\beta_{n-1})$, if $g_{k+1}\beta_kg_{k}^{-1} = \beta^0_k$.
By the definition of $\mm^{-1}_{\surj}(0)$, for any point of $\mm^{-1}_{\surj}(0)$, we can always find fitting matrices for this point.
But the choice of fitting matrices for a given point in $\mm^{-1}_{\surj}(0)$ is not unique.

Let $g_0,\cdots,g_n$ be the fitting matrices for $(\alpha_0,\beta_0,\cdots,\alpha_{n-1},\beta_{n-1})\in \mm^{-1}_{\surj}(0)$.
Using (\ref{eq:slmm}), by induction, we can show that $g_0\alpha_0\beta_0g_0^{-1}$ is of the following form, c.f.~\cite[(7.1)]{Dancer:2013aa}.
\begin{equation}
  \label{eq:defX}
  g_0\alpha_0\beta_0g_0^{-1} =
\begin{pmatrix}
    0      & *      &        & \cdots      & * \\
    0      & \rho_1    & \ddots     &         & \vdots \\
    \vdots      &  \ddots    & \ddots     &  \ddots     & \\
           &        &  \ddots    & \ddots    & * \\
    0      & \cdots      &        & 0      & \rho_n
  \end{pmatrix},
\end{equation}
where $\rho_k = \sum_{p=1}^{k}\tau_p$.
In general, $g_0\alpha_0\beta_0g_0^{-1}$ is not trace-free, consequently, not an element of $\kb$.
But we can use the following homomorphism of Lie algebras
\begin{equation*}
  \Gamma(Y) \coloneqq Y - \frac{\Tr{Y}}{n+1} I_{n+1}:\; \kgl(n+1,\bC)\rightarrow \ksl(n+1,\bC),
\end{equation*}
which is induced from the natural decomposition $\GL{n+1} = \SL{n+1}\times \bC^*$.
With such a homomorphism, we see that $\Gamma(g_0\alpha_0\beta_0g_0^{-1}) \in \kb$.

With the above settings, we define the following map,\footnote{If we use the notations of \S~\ref{sec:dks-type-varieties}, the orbit equivalence class $[(\alpha_0,\beta_0,\cdots,\alpha_{n-1},\beta_{n-1})]$ should be written as $\pi((\alpha_0,\beta_0,\cdots,\alpha_{n-1},\beta_{n-1}))$. In this section, we prefer using this shorter notation.}
\begin{equation}
  \label{eq:def-phi}
  \begin{aligned}
    \mm^{-1}_{\surj}(0)/\rSL_{\bfv} &\rightarrow G\times_U \kb \simeq \rT^*(G/U) \\
    [(\alpha_0,\beta_0,\cdots,\alpha_{n-1},\beta_{n-1})]  &\mapsto [g_0^{-1}, \Gamma(g_0\alpha_0\beta_0g_0^{-1})].
  \end{aligned}
\end{equation}
Note that to define the above map, for each point in $\mm^{-1}_{\surj}(0)$, we need to choose a family of fitting matrices for this point.
We will show that the image of the map doesn't depend on the different choices of fitting matrices.

Let $\aset{g'_i}$ be another family of matrices fitting for the point $(\alpha_0,\beta_0,\cdots,\alpha_{n-1},\allowbreak \beta_{n-1})$.
Using the definition of fitting matrices, one can check that there exists $u \in U$ such that $g'_0 = ug_0$.
Therefore, we have
\begin{equation*}
  ({g'}_0^{-1},g'_0\alpha_0\beta_0{g'}_0^{-1}) = u\cdot (g_0^{-1},g_0\alpha_0\beta_0g_0^{-1}),
\end{equation*}
which implies that the map (\ref{eq:def-phi}) is well-defined.

By~\cite[Theorem~7.18]{Dancer:2013aa}, the map (\ref{eq:def-phi}) can be extended to an isomorphism between two varieties in a suitable way, which justifies the following definition.
As a result, we have a quiver theoretic construction of the variety $\afbc$.
\begin{definition}
  \label{def:phi}
  We define 
  \begin{equation*}
    \Xi: \cM(\bfv,\bfw) = \mm^{-1}(0)\dsl \rSL_{\bfv} \simarrow \afbc
  \end{equation*}
  to be the unique isomorphism whose restriction on the open dense subvariety $\mm^{-1}_{\surj}(0)/\rSL_{\bfv} \subseteq \cM(\bfv,\bfw)$ coincides with the map (\ref{eq:def-phi}).
\end{definition}
For our application, we can always restrict $\Xi$ to the open dense subset $\mm^{-1}_{\surj}(0)/\rSL_{\bfv}$.
Namely, in the following discussion, we always use the map (\ref{eq:def-phi}) to calculate $\Xi$ explicitly.
Moreover, since for any regular and semi-simple element $y \in \kb$, there exists $u_y\in U$ such that $\Ad{u_y}{y}$ is a diagonal matrix, using (\ref{eq:slmm}), we can show that $\Xi^{-1}\big(\rT^*(G/U)_{\ras}\big) \subseteq \mm^{-1}_{\surj}(0)/\rSL_{\bfv}$.

\begin{remark*}
  Using a result by Grosshans~\cite{Grosshans_1997aa}, as well as by Popov and Vinberg~\cite{Vinberg_1972aa}, we know that $\afbc$ is a normal variety.
  Therefore, the above isomorphism $\Xi$ implies that $\cM(\bfv,\bfw)$ is normal in this special case.
  In view of the normality of quiver varieties,~\cite{Crawley-Boevey_2003aa}, it seems reasonable to expect that the normality also holds for other DKS type varieties.
\end{remark*}

\paragraph{The equivariance properties of $\Xi$.}
As we have recalled, there are some natural group actions on $\afbc$.
We will see that the same groups can also act on $\cM(\bfv,\bfw)$ and further show that $\Xi$ is equivariant with respect to these group actions.
Since $G = \SL{n+1}$, we choose $H$ to be the Cartan subgroup of $G$ consisting of diagonal matrices.

First, we recall the $G\times H$ acts on $G\times_U \kb$ (on $\afbc$ consequently) in the following way, for $g\in G$, $h\in H$, $[g',x]\in G\times_U \kb$,
\begin{equation*}
  (g,h) \cdot [g',x] \coloneqq [gg'h, \Ad{h^{-1}}{x}].
\end{equation*}
Via the isomorphism $\Theta_s$, the above group action is just the $G\times T$ action on $\rT^*\widetilde{\cB}$ (identifying $H$ with $T$).

The $G$ action on $\cM$ is induced from the standard $G$ action on $\bC^{n+1}$, that is, for $g\in G$, 
\begin{equation*}
  g \cdot [(\alpha_0,\beta_0,\cdots,\alpha_{n-1},\beta_{n-1})] \coloneqq [(g\alpha_0,\beta_0g^{-1},\cdots,\alpha_{n-1},\beta_{n-1})],
\end{equation*}
where $[(\alpha_0,\beta_0,\cdots,\alpha_{n-1},\beta_{n-1})] \in \mm^{-1}_{\surj}(0)/\rSL_{\bfv}$.
By (\ref{eq:def-phi}), we know that $\Xi$ is a $G$-equivariant isomorphism.

To define the $H$ action on $\cM$, we need the following fractional homomorphism between $H$ and $(\bC^*)^n$,
\begin{align*}
  \varphi:H & \rightarrow (\bC^*)^n \\
  h = \diag{z_0,\cdots,z_{n-1},z_n} & \mapsto \big(\prod^n_{k=1}z_k^{1/n},\prod^n_{k=2}z_k^{1/(n-1)}, \cdots,z_{n-1}^{1/2}z_n^{1/2},z_n\big).
\end{align*}
Let $c_k$, $1\le k \le n$, be a component of $\varphi(h)$.
The $H$ action on $\cM$ is defined to be
\begin{equation*}
  h\cdot [(\alpha_0,\beta_0,\cdots,\alpha_{n-1},\beta_{n-1})] \coloneqq [\alpha_0c_1,c_1^{-1}\beta_0,\cdots,c_{n-1}^{-1}\alpha_{n-1}c_n,c_{n}^{-1}\beta_{n-1}c_{n-1}],
\end{equation*}
where $[(\alpha_0,\beta_0,\cdots,\alpha_{n-1},\beta_{n-1})] \in \mm^{-1}_{\surj}(0)/\rSL_{\bfv}$.

Since there are some ambiguities in the definition of $\varphi$, we need to check that the above $H$ action is well-defined.
Note that each equivalence class in $\mm^{-1}_{\surj}(0)/\rSL_{\bfv}$ has a representative element of the form,
\begin{equation*}
  [(g\alpha_0,\beta^0_0g^{-1},\cdots,\alpha_{n-1},\beta^0_{n-1})],
\end{equation*}
and the $H$ action on such elements is easy to calculate.
Choose $(g_1,\cdots,g_n)\in \rSL_{\bfv}$ associated with $h$ as follows,
\begin{equation*}
  g_k \coloneqq c_k h_k^{-1},\; h_k \coloneqq \diag{z_k,\cdots,z_n}.
\end{equation*}
We have
\begin{equation}
  \label{eq:ha-m}
  \begin{aligned}
    & h \cdot [(g\alpha_0,\beta^0_0g^{-1},\cdots,\alpha_{n-1},\beta^0_{n-1})] \\
    =& [(g\alpha_0c_1g_1^{-1}, g_1c_1^{-1}\beta^0_0g^{-1}, \cdots, g_{n-1}c_{n-1}^{-1}\alpha_{n-1}c_ng_n^{-1}, g_nc_{n}^{-1}\beta^0_{n-1}c_{n-1}g_{n-1}^{-1})]\\
    =& [(g\alpha_0h_1,h_1^{-1}\beta^0_0g^{-1},\cdots, h_{n-1}^{-1}\alpha_{n-1}h_n, \beta^0_{n-1})]
\end{aligned}
\end{equation}
Note that the last formula in the above equality doesn't depend on the choice of $\varphi(h)$, that is, it has no ambiguity at all.
Hence, we know that the $H$ action on $\cM$ is well-defined.
In fact, we can define the $H$ action on $\cM$ by simply using the last formula of (\ref{eq:ha-m}).
However, the advantage of the definition using $\varphi$ is that we can check the following relation between the $H$ action and the Weyl group action on $\cM$ easily.
For $[B]\in \cM$, $w\in \bW\simeq W_H$ and $h\in H$,
\begin{equation*}
  w(h\cdot [B]) = w(h) \cdot w([B]).
\end{equation*}
The details are left to readers.

On the other hand, by (\ref{eq:def-phi}) and (\ref{eq:ha-m}), we have
\begin{align*}
  & h\cdot \Xi([(g\alpha_0,\beta^0_0g^{-1},\cdots,\alpha_{n-1},\beta^0_{n-1})]) \\
  = & h\cdot [g, \Gamma(\alpha_0\beta^0_0)] = [gh, \Gamma(h^{-1}\alpha_0\beta^0_0h)]\\
  = & \Xi([g\alpha_0h_1,h_1^{-1}\beta^0_0g^{-1},\cdots, h_{n-1}^{-1}\alpha_{n-1}h_n, \beta^0_{n-1}]) \\
  = & \Xi(h\cdot [(g\alpha_0,\beta^0_0g^{-1},\cdots,\alpha_{n-1},\beta^0_{n-1})]),
\end{align*}
that is, $\Xi$ is $H$-equivariant.

\begin{remark*}
  The careful readers may have noticed that here we only define the $G$ and $H$ action on an open dense subvariety of $\cM$ actually.
  However, the equivalence property of $\Xi$ implies that these two actions can be extended to $\cM$.
  It is also possible to show this fact directly by using similar arguments as in~\cite{Ginzburg:2018la,Dancer:2013aa}.
\end{remark*}

\paragraph{The dependence of $\Xi$ on base points.}
Recall that the isomorphism $\Theta_{\sfs}$ between $\rT^*(G^{\mathrm{ad}}/U)$ and $\rT^*\widetilde{\cB}$ depends on an auxiliary parameter $\sfs\in \bO(\kb)$.
The isomorphism $\Xi$ also depends on a parameter: the choice of base points, $\beta^0_k$, $0\le k \le n-1$.
Since $\bO(\kb)$ is a $T$-torsor, the set $\aset{\Theta_{\sfs}}$ is also a $T$-torsor.
Via varying the base points, we will see that a similar result is also true for $\Xi$.

Choose another family of matrices $\aset{\bar{\beta}^0_k}$ as the base point, such that $\bar{\beta}^0_k$ is of the form $(\zm{n-k}{1}\, |\, F_{n-k})$, where $F_{n-k}$ is an invertible upper triangular matrix of size $n-k$.
By using $\aset{\bar{\beta}^0_k}$ to replace $\aset{{\beta}^0_k}$ when defining the fitting matrices, using the same method as (\ref{eq:def-phi}), we can define another isomorphism $\bar{\Xi}$.

Moreover, we can find $a_k\in \bC^*$, $0 \le k \le n-1$ so that $\aset{(\zm{n-k}{1}\, |\, a_kI_{n-k})}$ and $\aset{\bar{\beta}^0_k}$ yield the same isomorphism $\bar{\Xi}$.
Therefore, the collection of all possible isomorphisms $\Xi$, constructed by using different base points, constitutes a $T$-torsor.

\begin{remark}[A notation convention]
  \label{rk:nc}
  In the remaining part of this paper, the symbol $\Xi$ is only reserved for the isomorphism given in Definition~\ref{def:phi}, that is, the isomorphism defined using the base points $\aset{\beta^0_k}$.
  Meanwhile, we also choose a special $\Theta_{\sfs}$ by choosing $\sfs$ to be the element such that $e_1 = \cdots = e_n = 1$, where $e_i$ is the symbol used in Example of \S~\ref{sub:goveru}.
  We will write this special $\Theta_{\sfs}$ as $\Theta$ for short.
  As we have explained in Remark~\ref{rk:sc}, for $G= \SL{n+1}$, $\Theta$ is only a finite covering map, with fibers isomorphic to $\bZ_{n+1}$.
\end{remark}

\subsection{The matrix representation of $S_k$.}

By identifying $\afbc$ with $\cM(\bfv,\bfw)$, we can calculate explicitly the Weyl group action of $\cM(\bfv,\bfw)$.
More precisely, for $[g,x]\in \rT^*(G/U)_{\ras}$, $B\in \mm^{-1}_{\surj}(0)$ and $\Xi([B]) = [g,x]$, we define
\begin{equation}
  \label{eq:sk-m}
  \sigma_k([g,x]) \coloneqq \Xi(S_k([B])).
\end{equation}
We are going to calculate $\sigma_k([g,x])$ using a specific matrix.

To make the calculation simpler, for $[g,x]\in \rT^*(G/U)_{\ras}$, we will always assume that $x$ is a diagonal matrix, which causes no loss of generality since we are dealing with regular and semi-simple elements.
Then we will show the following result.

\begin{proposition}
  \label{prop:skcal}
  For any $[g,x]\in \rT^*(G/U)$ such that $x = \diag{\lambda_0,\cdots,\lambda_n}$ is regular, we have $\sigma_k([g,x]) = [gW_k^{-1},\Ad{W_k}{x}]$, where
  \begin{equation}
    \label{eq:defwk}
    W_k = \mleft(
    \begin{array}{c|cc|c}
      I_{k-1}& \zm{k-1}{1} & \zm{k-1}{1} & \zm{k-1}{n-k} \\
      \hline
      \zm{1}{k-1} & 0 & (\lambda_k-\lambda_{k-1})^{-1} & \zm{1}{n-k} \\
      \zm{1}{k-1} & \lambda_{k-1}-\lambda_k & 0 & \zm{1}{n-k} \\
      \hline
      \zm{n-k}{k-1} & \zm{n-k}{1} & \zm{n-k}{1} & I_{n-k}
    \end{array}
    \mright)\,\in \SL{n+1}.
  \end{equation}
\end{proposition}

Before discussing the proof of this proposition, we comment that as an application of this explicit expression of $\sigma_k([g,x])$, one can check directly that the $\sigma_k$ induces a $\bW$ action on $\afbc$ without using Theorem~\ref{thm:main-re} and the isomorphism $\Xi$.

\begin{proof}
  During the proof, we will simply take $g = e$ because the proof for the general case is essentially same.
  As a preparation, we try to find $B = (\alpha_0,\beta_0,\cdots,\alpha_{n-1},\allowbreak \beta_{n-1}) \in \mm^{-1}_{\surj}(0)$ satisfying $\Xi([B]) = [e,x]$.
  
  Since $x = \diag{\lambda_0,\cdots,\lambda_n}$, inspired by (\ref{eq:defX}), we choose $X \coloneqq \diag{0,\rho_1,\cdots,\allowbreak\rho_n}$ such that $x = \Gamma(X)$ and $\tau_p$, $1\le p \le n$, to be the complex numbers satisfying $\rho_l = \sum_{p=1}^{l}\tau_p$.
  Then, by setting $\rho_0 =0$, we have
  \begin{equation}
    \label{eq:trl}
    \tau_p = \rho_p- \rho_{p-1} = \lambda_{p} - \lambda_{p-1},\; 1 \le p \le n.
  \end{equation}
  
  Moreover, by our construction of $\Xi$, we can choose $\beta_l =\beta^0_l = (\zm{n-l}{1}\, |\, I_{n-l})$ and $\alpha_0$ must satisfy the following equality,
  \begin{equation}
    \label{eq:xab}
    X = \alpha_0\beta_0.
  \end{equation}
  Then, combining (\ref{eq:slmm}) and (\ref{eq:xab}), we know that by taking
  \begin{equation}
    \label{eq:alphak}
    \alpha_l = \mleft(
    \begin{array}{c}
      \zm{1}{n-l} \\
      \hline
      C_{n-l}
    \end{array}\mright),\quad
    C_{n-l} \coloneqq \diag{\tau_{l+1},\tau_{l+1}+ \tau_{l+2},\cdots,\textstyle\sum_{p=l+1}^n\tau_p},
  \end{equation}
  we have $\Xi([B]) = [e,x]$.

  From now on, we fix a vertex $k$, $1\le k \le n$ and use the notations in \S~\ref{sub:weyl-group-action}.
  In the current situation, $T_k$ is equal to $\bC^{n+2-k} \oplus \bC^{n-k}$ and the maps $\om_k(B),\im_k(B)$ are written as
  \begin{equation*}
   \om_k(B) =
    \begin{pmatrix}
      \alpha_{k-1} \\ \beta_k
    \end{pmatrix},\quad
    \im_k(B) =
    \begin{pmatrix}
      \beta_{k-1} & -\alpha_k
    \end{pmatrix}.
  \end{equation*}

  To calculate $S_k([B])$, we need to find $B'\in \mm^{-1}(0)$ such that $(B,B')\in Z_k$.
  In fact, by (\ref{eq:alphak}), one can verify that if $\om_k(B')$ takes the following value, $(B,B')$ will satisfy the condition~\ref{con:c2} in Definition~\ref{def:zk},
  \begin{equation}
    \label{eq:aprime}
    \om_k(B') =
    \begin{pmatrix}
      \alpha'_{k-1} \\ \beta'_k
    \end{pmatrix},\text{ where }
    \alpha'_{k-1} = \mleft(
    \begin{array}{cc}
      1 & \zm{1}{n-k} \\
      0 & \zm{1}{n-k} \\
      \zm{n-k}{1} & C_{n-k}
    \end{array}
    \mright)\,,\;
\beta_k' = \beta_k.
  \end{equation}
  Choosing $\om_k(B')$ as above, if $(B,B')$ further satisfies the condition~\ref{con:c3} in Definition~\ref{def:zk}, $\im_k(B') = (\beta'_{k-1},-\alpha'_k)$ must take the value
  \begin{equation}
    \label{eq:bprime}
    \beta'_{k-1} = \mleft(
    \begin{array}{ccc}
      -\tau_k & 0 & \zm{1}{n-k} \\
      \zm{n-k}{1} & \zm{n-k}{1} & I_{n-k}
    \end{array}
    \mright),\quad
    \alpha'_{k} =
    \begin{pmatrix}
      \zm{1}{n-k} \\
      C_{n-k} + \tau_kI_{n-k}
    \end{pmatrix}.
  \end{equation}
  Due to (\ref{eq:aprime}) and (\ref{eq:bprime}), by choosing
  \begin{equation*}
    B' = (\alpha_0,\beta_0,\cdots, \alpha'_{k-1},\beta'_{k-1}, \alpha'_k,\beta'_k,\cdots,\alpha_{n-1},\beta_{n-1}),
  \end{equation*}
  we have $(B, B')\in Z_k$.
  Therefore, by (\ref{eq:sk-m}), the following equality holds,
  \begin{equation}
    \label{eq:sk-m-1}
    \sigma_k([e,x]) = \Xi(S_k([B])) = \Xi([B']).
  \end{equation}

  The next step is to calculation $\Xi([B'])$.
  Take $g_p$, $0\le p \le n$, as follows,
\begin{gather*}
    g_{k-1} = \mleft(
    \begin{array}{cc|c}
0 & \tau_k^{-1} & \zm{1}{n-k}\\
      -\tau_k & 0 & \zm{1}{n-k} \\
      \hline
      \zm{n-k}{1} & \zm{n-k}{1} & I_{n-k}
    \end{array}
    \mright)\,, \\
g_p =
    \begin{pmatrix}
      1 & \zm{1}{n-p} \\
      \zm{n-p}{1} & g_{p+1} 
    \end{pmatrix},\text{ for }0\le p \le k-2, \text{ and }
    g_p = I_{n-p+1}\text{ for }k\le p \le n.
  \end{gather*}
Due to (\ref{eq:defwk}) and (\ref{eq:trl}), we have $W_k = g_0$.
  With the above definition of $g_p$, we can check that, for $0\le p \le n-1$,
  \begin{equation}
    \label{eq:gp}
    \begin{gathered}
      g_{p+1}\beta_pg_{p}^{-1} = \beta^0_p\;\text{ if }p\notin \aset{k-1, k},\\
      g_{p+1}\beta'_pg_{p}^{-1} = \beta^0_p\;\text{ if }p\in \aset{k-1, k}.
    \end{gathered}
  \end{equation}
  
By (\ref{eq:gp}), we can use (\ref{eq:def-phi}) to calculate $\Xi([B'])$.
  We need to discuss two cases depending on whether $k = 1$ or not.

  \medskip
  \noindent \uline{\textsc{Case 1:} $2\le k \le n$.}
  In this case, the first two components of $B'$ are $\alpha_0$ and $\beta_0$.
  Therefore, by (\ref{eq:def-phi}), (\ref{eq:xab}) and (\ref{eq:gp}), we have
  \begin{multline}
    \label{eq:sk-m-2}
    \Xi([B']) = [W^{-1}_k, \Gamma(\Ad{W_k}{\alpha_0\beta_0})] = [W^{-1}_k,\Ad{W_k}{\Gamma(\alpha_0\beta_0)}] \\
    = [W^{-1}_k,\Ad{W_k}{\Gamma(X)}] = [W^{-1}_k, \Ad{W_k}{x}].
  \end{multline}

  \medskip
  \noindent \uline{\textsc{Case 2:} $k = 1$.}
  In this case, the first two components of $B'$ are $\alpha'_0$ and $\beta'_0$, which are given in (\ref{eq:aprime}), (\ref{eq:bprime}).
  Hence, we have
  \begin{equation}
    \label{eq:sk-m-3}
    \alpha'_0\beta'_0 =\mleft(
    \begin{array}{cc|c}
      -\tau_1 & 0 & \zm{1}{n-1} \\
      0 & 0 & \zm{1}{n-1} \\
      \hline
      \zm{n-1}{1} & \zm{n-1}{1} & C_{n-1}
    \end{array}
    \mright)\,= X - \tau_1I_{n+1}.
  \end{equation}
  Due to the above equality, by (\ref{eq:def-phi}) and (\ref{eq:gp}), we know that
  \begin{multline*}
    \Xi([B']) = [W^{-1}_1, \Gamma(\Ad{W_1}{\alpha'_0\beta'_0})] = [W^{-1}_1,\Gamma(\Ad{W_1}{X} - \tau_1I_{n+1})] \\
    = [W^{-1}_1,\Ad{W_1}{\Gamma(X)}] = [W^{-1}_1, \Ad{W_1}{x}].
  \end{multline*}

  By (\ref{eq:sk-m-1}), (\ref{eq:sk-m-2}) and (\ref{eq:sk-m-3}), the proof of Proposition~\ref{prop:skcal} finishes.

\end{proof}
\subsection{The proof of Theorem~\ref{thm:main2}.}
\label{sub:pf-thm2}
Note that the map $\Theta$ induces an isomorphism between $\big(\rT^*(G^{\mathrm{ad}}/U)\big)_{\aff}$ and $(\rT^*\widetilde{\cB})_{\aff}$.
Therefore, the \grk\ action induces a $\bW$ action on $\big(\rT^*(G^{\mathrm{ad}}/U)\big)_{\aff}$.
On the other hand, via the isomorphism $\Xi$, the Weyl group action on $\cM$ can also induce a group action on $\afbc$.
In this subsection, we compare these two Weyl group actions, that is, we will prove Theorem~\ref{thm:main2}.

As before, for $k\in I$, let $s_k\in \bW$ be a simple reflection.
Recall that for any $p\in \rT^*\widetilde{\cB}_{\ras}$, $s_k(p)$ has been defined in \S~\ref{sub:weyl-gu} in an explicit way.
We restate Theorem~\ref{thm:main2} in the following form.
\begin{theorem}[{$=$Theorem~\ref{thm:main2}}]
  \label{thm:main2-re}
  For any $[g,x]\in \rT^*(G/U)$ such that $x = \diag{\lambda_0,\allowbreak \cdots,\allowbreak \lambda_n}$ is regular, we have
  \begin{equation*}
    s_k(\Theta([g,x])) = \Theta([gW_k^{-1},\Ad{W_k}{x}]) = \Theta(\sigma_k([g,x])).
  \end{equation*}
  Consequently, the map $\Theta \circ \Xi: \cM(\bfv,\bfw) \rightarrow (\rT^*\widetilde{\cB})_{\aff}$ is $\bW$-equivariant.
\end{theorem}
\begin{remark*}
  In fact, it is easy to further verify that $\Theta \circ \Xi$ is $\bW \ltimes T$-equivariant.
\end{remark*}

\begin{proof}
  Due to Proposition~\ref{prop:skcal} and the definition of $\Theta$, to prove Theorem~\ref{thm:main2-re}, we only need to show
  \begin{equation}
    \label{eq:weyleq1}
    s_k\big(\Ad{g}{\kb},\Ad{g}{\sfs},\Ad{g}{x}\big) = \Theta([gW_k^{-1},\Ad{W_k}{x}]).
\end{equation}
  As before, it suffices to show the above equality with the assumption $g= e$.

  To begin with, we notice that since $x$ is regular, we can find $u\in U$ such that
  \begin{equation}
    \label{eq:yux}
    y = \Ad{u}{x} =
\begin{pmatrix}
      \lambda_0    & 1      & 0      & \cdots      &        & 0 \\
      0      & \ddots     &   \ddots    & \ddots     &        & \vdots \\
      \vdots      & \ddots     &   \ddots    & \ddots     &  \ddots     & \\
      &        &  \ddots     &   \ddots   &  \ddots     & 0 \\
      &        &        &   \ddots    &   \ddots     & 1 \\
      0      &        & \cdots      &        &    0   & \lambda_n
    \end{pmatrix}.
  \end{equation}
  Let $\bkb$ be the Borel subalgebra of the lower triangular matrices.
  Thus, $\kb,\bkb$ are in opposite position.
  Since we have chosen $\sfs$ to be of a special form, c.f.\ Remark~\ref{rk:nc}, by (\ref{eq:yux}), we have
  \begin{equation}
    \label{eq:ys}
    \kappa_{\bkb,\kb}(y \bmod \bkb) = \sfs.
  \end{equation}
  Therefore, by (\ref{eq:yux}) and (\ref{eq:ys}), we get
  \begin{multline}
    \label{eq:xs}
    \kappa_{\ukb,\kb}(x \bmod \ukb) = \kappa_{\ukb,\kb}(\Ad{u^{-1}}{y} \bmod \ukb)\\
    = \Ad{u^{-1}}{\sfs} = \sfs.
  \end{multline}
  Combining (\ref{eq:xs}) and the definition of $\kappa$ in Proposition~\ref{prop:gk}, we know that
  \begin{equation}
    \label{eq:bsx}
(\kb, \sfs, x) = \kappa(\ukb, \kb, x).
  \end{equation}
  Hence, by (\ref{eq:defw1}) and (\ref{eq:bsx}), we have
  \begin{equation}
    \label{eq:weyleq2}
    s_k({\kb},{\sfs},{x}) = \kappa(\ukb, \kb^{s_k^{-1}}, x).
  \end{equation}

  To calculate the r.h.s.\ of (\ref{eq:weyleq2}), we use the following matrix,
  \begin{equation*}
    r \coloneqq u^{-1}
    \begin{pmatrix}
      0      & 1    & 0      & \cdots       & 0 \\
      \vdots      & \ddots      &    \ddots  &   \ddots   & \vdots \\
             &         &  \ddots     &   \ddots    & 0 \\
      \vdots      &        &        &  \ddots      & 1 \\
      0      &        & \cdots      &         & 0
    \end{pmatrix}u \in \kd\big(\Ad{u^{-1}}{\bkb}\big) \cap \ku(\kb).
  \end{equation*}
  Due to (\ref{eq:yux}) and (\ref{eq:ys}), one can check that
  \begin{subequations}
    \begin{align}
      r \bmod \Ad{u^{-1}}{\bkb} &= x \bmod \Ad{u^{-1}}{\bkb}, \label{eq:r1}\\
      r \bmod \kb &= \sfs.\label{eq:r2}
    \end{align}
  \end{subequations}

Now, we use Lemma~\ref{lm:nk} to find $n_k \in\bar{U} = [\bar{B},\bar{B}]$ such that $W_ku^{-1}n_ku \in U$.
  By the definition of $\kb^{s_k^{-1}}$ in \S~\ref{sub:weyl-gu}, we have $\kb^{s_k^{-1}} = \Ad{W_k^{-1}}{\kb}$, which implies
  \begin{equation}
    \label{eq:kbs}
    \kb^{s_k^{-1}} = \Ad{u^{-1}n_ku}{\kb}.
  \end{equation}
Since $n_k\in \bar{U}$, by (\ref{eq:kbs}) and the definition of $r$, we have
  \begin{subequations}
    \begin{gather}
    \Ad{u^{-1}n_ku}{r} \in \kd\big(\Ad{u^{-1}}{\bkb}\big) \cap \ku(\kb^{s_k^{-1}}),\label{eq:adr}\\
    \Ad{u}{r} - \Ad{n_ku}{r} \in \bkb.\label{eq:adr2}
  \end{gather}
  \end{subequations}

  By (\ref{eq:r1}) and (\ref{eq:adr2}), one can check that the following equality holds,
  \begin{equation}
    \label{eq:use-r}
    x \bmod \Ad{u^{-1}}{\bkb} = r \bmod \Ad{u^{-1}}{\bkb} {=} \Ad{u^{-1}n_ku}{r} \bmod \Ad{u^{-1}}{\bkb}.
  \end{equation}
On the other hand, by (\ref{eq:r2}) and (\ref{eq:kbs}), we have
\begin{equation}
    \label{eq:use-r1}
    \Ad{u^{-1}n_ku}{\sfs} = \Ad{u^{-1}n_ku}{r} \bmod \Ad{u^{-1}n_ku}{\kb} = \Ad{u^{-1}n_ku}{r} \bmod \kb^{s_k^{-1}}.
  \end{equation}
  Combining (\ref{eq:adr}), (\ref{eq:use-r}) and (\ref{eq:use-r1}), we have
  \begin{equation}
    \label{eq:kbb}
    \kappa_{\ukb,\kb^{s_k^{-1}}}(x \bmod \Ad{u^{-1}}{\bkb}) = \Ad{u^{-1}n_ku}{r} \bmod \kb^{s_k^{-1}} = \Ad{u^{-1}n_ku}{\sfs}.
  \end{equation}

  By (\ref{eq:kbs}), (\ref{eq:kbb}) and the definition of $\kappa$ in Proposition~\ref{prop:gk}, one concludes that
  \begin{equation}
    \label{eq:bsx2}
    \kappa(\ukb,\kb^{s_k^{-1}}, x) = (\kb^{s_k^{-1}}, \Ad{u^{-1}n_ku}{\sfs}, x) = \Theta([u^{-1}n_ku,\Ad{u^{-1}n_k^{-1}u}{x}]).
  \end{equation}
  By (\ref{eq:weyleq2}), (\ref{eq:bsx2}) and $W_ku^{-1}n_ku \in U$, we have shown (\ref{eq:weyleq1}) and the proof of Theorem~\ref{thm:main2-re} finishes.
\end{proof}

During the above proof, we have used the following lemma.
\begin{lemma}
  \label{lm:nk}
  Using the same notations in the proof of Theorem~\ref{thm:main2-re}, there exists $n_k \in\bar{U} = [\bar{B},\bar{B}]$ such that $W_ku^{-1}n_ku \in U$.
\end{lemma}
\begin{proof}
Before running into details of the proof, we would like to remind readers that a weaker version of Lemma~\ref{lm:nk}, that is, there exists $n_k\in \bar{U}$ such that $W_ku^{-1}n_ku \in B$, is easy to verify.
We can deduce it from the following fact: the set $\Omega$ consisting of pairs of Borel subalgebras in opposite position is the unique open dense orbit of the diagonal action of $G$ on $\cB \times \cB$.

  To show the little stronger result claimed in the lemma, we need some calculation.
  In fact, we will try to find $n_k\in \bar{U}$ such that $n_k^{-1}uW_k^{-1}u^{-1}\in U$.
  Using (\ref{eq:defwk}), $W_k^{-1}$ can be written in the following block form,
  \begin{equation*}
    W_k^{-1}=
    \begin{pmatrix}
      I_{k-1} &  & \\
      & F & \\
      & & I_{n-k}
    \end{pmatrix},\quad\text{where}\quad
    F \coloneqq
    \begin{pmatrix}
      0 & (\lambda_{k-1}-\lambda_{k})^{-1}\\
      \lambda_k-\lambda_{k-1} & 0 
    \end{pmatrix}.
  \end{equation*}
  Using (\ref{eq:yux}), we can also write $u$ as a block matrix with the same type of $W_k^{-1}$,
  \begin{equation*}
    u =
    \begin{pmatrix}
      A_{11} & * & * \\
      & A_{22} & * \\
      & & A_{33}
    \end{pmatrix},\quad\text{where}\quad
    A_{22} \coloneqq
    \begin{pmatrix}
      1 & (\lambda_k-\lambda_{k-1})^{-1} \\
      0 & 1
    \end{pmatrix}.
  \end{equation*}
  Similarly, let $n_k^{-1}$ be the following block matrix, again, with the same type of $W_k^{-1}$,
  \begin{equation*}
    n_k^{-1} \coloneqq
    \begin{pmatrix}
      I_{k-1} &  & \\
      & N & \\
      &  & I_{n-k}
    \end{pmatrix},\quad\text{where}\quad
    N \coloneqq
    \begin{pmatrix}
      1 & 0 \\
      \lambda_{k-1}-\lambda_k & 1
    \end{pmatrix}.
  \end{equation*}

  By this definition, $n_k^{-1}$ (or $n_k$) lies in $\bar{U}$.
  We assert that $n_k^{-1}uW_k^{-1}u^{-1}\in U$, i.e.\ $n_k$ satisfying the requirement of Lemma~\ref{lm:nk}.
  In fact, by using the block matrix expressions for $n_k^{-1}$, $u$ and $W_k^{-1}$ as above, we have
  \begin{equation*}
    n_k^{-1}uW_k^{-1}u^{-1} = n_k^{-1} \cdot
    \begin{pmatrix}
      I_{k-1} & * & * \\
      & A_{22}FA_{22}^{-1} & * \\
      &  & I_{n-k}
    \end{pmatrix} =
    \begin{pmatrix}
      I_{k-1} & * & * \\
      & NA_{22}FA_{22}^{-1} & * \\
      &  & I_{n-k}
    \end{pmatrix},
  \end{equation*}
  where
  \begin{multline*}
    NA_{22}FA_{22}^{-1} =
    \begin{pmatrix}
      1 & 0 \\
      \lambda_{k-1}-\lambda_k & 1
    \end{pmatrix}
    \begin{pmatrix}
      1 & 2(\lambda_{k-1}-\lambda_{k})^{-1} \\
      \lambda_k-\lambda_{k-1} & -1
    \end{pmatrix} \\
    =
    \begin{pmatrix}
      1 & 2(\lambda_{k-1}-\lambda_{k})^{-1} \\
      0 & 1
    \end{pmatrix}.
  \end{multline*}
  Therefore, $n_k^{-1}uW_k^{-1}u^{-1}\in U$.
\end{proof}
\begin{remark}[The simply connected group case]
  \label{rk:pf-sc}
  By using $\widetilde{\cB}^{\mathrm{sc}}$ mentioned in Remark~\ref{rk:sc} to replace $\widetilde{\cB}$,
we can actually show that there is a $\bW$-equivariant isomorphism between $\cM(\bfv,\bfw)$ and $(\rT^*\widetilde{\cB}^{\mathrm{sc}})_{\aff}$.
  Although the method to show this stronger result is the same with the proof of Theorem~\ref{thm:main2-re}, on a technical level, some differences between these two proofs are unavoidable.
  In this remark, we collect the gists of the necessary modifications for the proof of this stronger result.

  As in Remark~\ref{rk:sc}, we choose $\rv_i$ to be a highest weight vector (with respect to $\kb$) for $V_{\varpi_i}$, the irreducible representation of $\SL{n+1}$ with the highest weight $\varpi_i$.
  To replace $\bO(\kb)$, we use $\bO^{\rsc}(\kb) \coloneqq \prod_i \bC^* \cdot \rv_i$.
  Accordingly, we replace $\bO_-(\bkb)$ by $\bO^{\rsc}_{-}(\bkb) \coloneqq \prod_i ((V_{\varpi_i}/\ku(\bkb)V_{\varpi_i}) - \aset{0})$.
  Using $\bO^{\rsc}(\kb)$ and $\bO^{\rsc}_{-}(\bkb)$, we can define $\widetilde{\cB}^{\rsc}$, $\widetilde{\cB}^{\rsc}_{-}$ and $\cZ^{\rsc}_{\ras}$ as before.
  Moreover, like $\kappa_{\bkb,b}$, we can define an isomorphism $\kappa^{\rsc}_{\bkb, \kb}: \bO_{-}(\bkb) \rightarrow \bO^{\rsc}(\kb)$, which is induced by the following maps,
  \begin{equation*}
    \bC \cdot \rv_i \hookrightarrow V_{\varpi_i} \twoheadrightarrow V_{\varpi_i}/\ku(\bkb)V_{\varpi_i}.
  \end{equation*}
  By fixing some isomorphisms as~\cite{Ginzburg:2018la}, there is a covering map $\pi_+: \bO^{\rsc}(\kb) \rightarrow \bO(\kb)$ (resp.\ $\pi_-:\bO^{\rsc}_{-}(\bkb) \rightarrow \bO_{-}(\bkb)$).
  One can check that these two covering maps make the following diagram commutes.
  \begin{equation*}
    \begin{tikzcd}
      \bO^{\rsc}_{-}(\bkb) \arrow[d, "\pi_-"] \arrow[r,"\kappa^{\rsc}_{\bkb,\kb}"] & \bO^{\rsc}(\kb) \arrow[d, "\pi_+"] \\
      \bO_{-}(\bkb) \arrow[r,"\kappa_{\bkb,\kb}"] & \bO(\kb).
    \end{tikzcd}
  \end{equation*}
  Like $\Theta$, to construct the isomorphism $\Theta^{\rsc}$ between $\afbc$ and $(\rT^*\widetilde{\cB}^{\mathrm{sc}})_{\aff}$, we need to choose a base point $\sfs^{\rsc}\in \bO^{\rsc}(\kb)$.
  Without loss of generality, we can choose $\sfs^{\rsc}$ to be $(\rv_1,\cdots,\rv_n)$ and assume that $\pi_+(\sfs^{\rsc}) = \sfs$.
  Using these notations, what we are going to show is that the isomorphism
  \begin{equation*}
    \Theta^{\rsc} \circ \Xi: \cM(\bfv,\bfw) \rightarrow (\rT^*\widetilde{\cB}^{\mathrm{sc}})_{\aff}
  \end{equation*}
  is $\bW$-equivariant.

  Let $\sfs^{\rsc}_x$ be the preimage of $\sfs^{\rsc}$ under $\kappa^{\rsc}_{\bkb,\kb}$.
  By definition,
  \begin{equation*}
    \sfs^{\rsc}_x = (\rv_1 \bmod \ku(\bkb)V_{\varpi_1}, \cdots, \rv_n \bmod \ku(\bkb)V_{\varpi_n}).
  \end{equation*}
  Then one can check that $\pi_-(\sfs^{\rsc}_x) = \Ad{u}{x} \bmod \bkb$.
  Using the notations in the proof of Theorem~\ref{thm:main2-re}, to prove the $\bW$-equivariance of $\Theta^{\rsc} \circ \Xi$, the crux is to show
  \begin{equation}
    \label{eq:use-r2}
    \kappa^{\rsc}_{\ukb,\kb^{s_k^{-1}}}(u^{-1}(\sfs^{\rsc}_x)) = (u^{-1}n_ku)(\sfs^{\rsc}),
  \end{equation}
  which corresponds to (\ref{eq:kbb}).
  By the definition of $\kappa^{\rsc}$ and $\sfs^{\rsc}$, (\ref{eq:use-r2}) is equivalent to
  \begin{equation*}
    n_k(\rv_i) - \rv_i \in \ku(\bkb)V_{\varpi_i}, \;1\le i \le n,
  \end{equation*}
  which is a counterpart of (\ref{eq:adr2}) and follows from the definition of $n_k$ in Lemma~\ref{lm:nk}.

\end{remark}

\bibliographystyle{amsplain}

\end{document}